
\magnification\magstep1  



\font\tenbbb=msbm10 \font\sevenbbb=msbm7 
\skewchar\tenbbb=127 \skewchar\sevenbbb=127 \newfam\bbbfam
\textfont\bbbfam=\tenbbb \scriptfont\bbbfam=\sevenbbb
\def\Bbb{\fam\bbbfam}

\font\eightrm=cmr8	


\def\cite#1{{\rm[#1]}}	
\def\eq#1{{\rm(#1)}}	
\def\opname#1{\mathop{\rm#1}\nolimits} 
\def\refno#1. #2\par{\smallskip\item{\rm #1.}#2\par} 
\def\today{\the\day\space \ifcase\month\or
January\or February\or March\or April\or May\or June\or
July\or August\or September\or October\or November\or December\fi
\space\the\year}	


\newtoks\secnum
\outer\def\section#1. #2\par{
\vskip 0pt plus.2\vsize	
\penalty -250 \vskip 0pt plus-.2\vsize
\bigskip\bigskip \vskip\parskip
\global\secnum={#1}\message{#1. #2}\leftline{\bf#1. #2}\nobreak
\smallskip\noindent}

\def\subsection#1.#2. #3\par{\medskip 
\noindent{\bf#1.#2. \it#3}\par\nobreak
\smallskip\noindent\ignorespaces}

\def\declare#1. #2\par{\medskip\noindent{\bf#1.}\rm\enspace
\ignorespaces#2\par\smallskip} 

\def\Proof.{\noindent{\it Proof}.\enspace\ignorespaces}
\def\endproof{\qed\par\smallskip} 

\def\qed{\allowbreak\qquad\null\nobreak\hfill\square}

\newtoks\leftrunhead \newtoks\rightrunhead

\headline={\ifnum\pageno>1 \ifodd\pageno \hfill
\eightrm \the\rightrunhead \hfill \llap{\tenrm\folio}\else
\rlap{\tenrm\folio}\hfill \eightrm \the\leftrunhead \hfill \fi
\else\hfil \fi}
\footline={\hfil} 
\hsize=16truecm 
\voffset=1truepc 


\newbox\ncintdbox \newbox\ncinttbox 
\setbox0=\hbox{$-$}
\setbox2=\hbox{$\displaystyle\int$}
\setbox\ncintdbox=\hbox{\rlap{\hbox
to \wd2{\hskip-.125em\box2\relax\hfil}}\box0\kern.1em}
\setbox0=\hbox{$\vcenter{\hrule width 4pt}$}
\setbox2=\hbox{$\textstyle\int$}
\setbox\ncinttbox=\hbox{\rlap{\hbox
to \wd2{\hskip-.175em\box2\relax\hfil}}\box0\kern.1em}

\let\Lpolish=\L	
\def\square{\vbox{\hrule\hbox{\vrule height5.2pt \hskip 5.2pt
\vrule}\hrule}}	
\def\stroke{\mathbin|}	



\def\a{\alpha}
\def\b{\beta}
\def\dl{\delta}	
\def\eps{\varepsilon}
\def\Ga{\Gamma}
\def\ga{\gamma}
\def\La{\Lambda}
\def\la{\lambda}
\def\Om{\Omega}
\def\om{\omega}


\def\D{{\cal D}}	
\def\H{{\cal H}}	
\def\K{{\cal K}}	
\def\L{{\cal L}}	
\def\M{{\cal M}}	
\def\Oh{{\cal O}}	
\def\SS{{\cal S}}	
\def\X{{\cal X}}	


\def\C{{\Bbb C}}	
\def\N{{\Bbb N}}	
\def\R{{\Bbb R}}	
\def\Sf{{\Bbb S}}	
\def\T{{\Bbb T}}	
\def\Z{{\Bbb Z}}	


\def\Dtr{\opname{Dtr}}	
\def\Fp{\opname{F.p.}}	
\def\Pf{\opname{Pf}}	
\def\rk{\opname{rk}}	
\def\Res{\opname{Res}}	
\def\spec{\opname{sp}}	
\def\supp{\opname{supp}}	
\def\Tr{\opname{Tr}}	
\def\vol{\opname{vol}}	
\def\Wres{\opname{Wres}}	
\def\wres{\opname{wres}}	


\def\del{\partial}	
\def\x{\times}	
\def\7{\dagger}	
\def\.{\cdot}	
\def\:{\colon}	



\def\shalf{{\scriptstyle{1\over2}}} 
\def\thalf{{\textstyle{1\over2}}}   
\def\tihalf{{\textstyle{i\over2}}}  


\def\textbf#1{{\bf#1}}	
\def\textit#1{{\it#1\/}}	
\def\textrm#1{{\rm#1}}	

\def\as{\qquad\hbox{as}\enspace} 
\def\sepword#1{\qquad\hbox{#1}\quad} 

\def\ala#1{\ifinner\mkern9mu \else\mkern18mu \fi(#1)} 
\def\fpart#1{\{#1\}}	
\def\piso#1{\lfloor#1\rfloor}	
\def\set#1{\{\,#1\,\}}	

\def\braket#1#2{(#1\stroke#2)} 
\def\ketbra#1#2{|#1)(#2|}	
\def\<#1,#2>{\langle#1,#2\rangle} 

\def\frac#1#2{{#1\over#2}}	
\def\pd#1#2{{\partial#1\over\partial#2}} 
\def\tfrac#1#2{{\textstyle{#1\over#2}}} 


\hyphenation{geo-me-tric geo-me-try ope-ra-tor ope-ra-tors pro-duct
pro-ducts pseudo-dif-fe-ren-tial skew-sym-met-ric sum-ma-bi-li-ty}


\def\Ackermann{1}
\def\Persephone{2}
\def\BransonG{3}
\def\Brownell{4}
\def\Carleman{5}
\def\CarminatiIKS{6}
\def\ChamConnes{7}
\def\ConnesA{8}
\def\ConnesNCGR{9}
\def\ConnesGrav{10}
\def\ElizaldeVZ{11}
\def\EstradaC{12}
\def\EstradaF{13}
\def\EstradaKfp{14}
\def\EstradaKbook{15}
\def\Amalthea{16}
\def\FullingA{17}
\def\GelfandShilov{18}
\def\Gilkey{19}
\def\GrossmannLS{20}
\def\Grubb{21}
\def\Gurarie{22}
\def\Hardy{23}
\def\HormanderS{24}
\def\HormanderIII{25}
\def\IochumKS{26}
\def\KalauW{27}
\def\KastlerEH{28}
\def\Lagrange{29}
\def\Gianni{30}
\def\Lojasiewicz{31}
\def\Cordelia{32}
\def\McKeanS{33}
\def\Mulholland{34}
\def\Rosenberg{35}
\def\Solomyak{36}
\def\Sirius{37}
\def\Wodzicki{38}


\rightline{\rm CPT--97/P.3452}
\rightline{\rm DFTUZ/97/04}
\rightline{\rm UCR--FM--9--97}
\rightline{\rm funct-an/9702001}

\bigskip\bigskip

\centerline{\bf ON SUMMABILITY OF DISTRIBUTIONS}
\medskip
\centerline{\bf AND SPECTRAL GEOMETRY}
\bigskip

\centerline{R. Estrada,$^1$ J. M. Gracia-Bond\'{\i}a$^{2\star}$
and J. C. V\'arilly$^{3\star}$}

\medskip

\centerline{\it $^1$P. O. Box 276, Tres R\'{\i}os, Costa Rica}
\smallskip
\centerline{\it $^2$Departamento de F\'{\i}sica Te\'orica,
Universidad de Zaragoza, 50009 Zaragoza, Spain}
\smallskip
\centerline{\it $^3$Centre de Physique Th\'eorique, CNRS--Luminy,
Case 907, 13288 Marseille, France}

\vfootnote{$^\star$}{On leave from:
Department of Mathematics, Universidad de Costa Rica,
2060 San Pedro, Costa Rica.\endgraf
Email: {\tt restrada@cariari.ucr.ac.cr}, {\tt varilly@cpt.univ-mrs.fr}}

\bigskip\medskip

\noindent{\bf Abstract.}\enspace
Modulo the moment asymptotic expansion, the Ces\`aro and parametric
behaviours of distributions at infinity are equivalent. On the
strength of this result, we construct the asymptotic analysis for
spectral densities, arising from elliptic pseudodifferential
operators. We show how Ces\`aro developments lead to efficient
calculations of the expansion coefficients of counting number
functionals and Green functions. The bosonic action functional
proposed by Chamseddine and Connes can more generally be validated as
a Ces\`aro asymptotic development.

\bigskip

\rightrunhead={ON SUMMABILITY OF DISTRIBUTIONS
AND SPECTRAL GEOMETRY} 

\leftrunhead={ESTRADA, GRACIA-BOND\'IA and V\'ARILLY}

\section 1. Introduction

Most approaches to spectral geometry rely on the asymptotic expansion
of the heat kernel and Tauberian theorems. In this work, motivated by
a string of recent papers by Connes, we develop spectral geometry from
a more fundamental object. According to a deep statement by
Connes~\cite{\ConnesGrav}, there is a one-to-one correspondence
between Riemannian spin geometries and commutative real $K$-cycles,
the dynamics of the latter being governed by the spectral properties
of its defining Dirac operator. On ordinary manifolds, gravity (of the
Einstein and the Weyl variety) is the only interaction naturally
described by the $K$-cycle~\cite{\Ackermann, \KalauW, \KastlerEH}.

That is to say, in noncommutative geometry, existence of gauge fields
requires the presence of a noncommutative manifold structure, whose
``diffeomorphisms'' incorporate the gauge transformations. Connes' new
gauge principle points thus to an intrinsic coupling between gravity
and the other fundamental interactions. In a recent
formulation~\cite{\ChamConnes}, the Yang--Mills action functional is
replaced by a ``universal'' bosonic functional of the form:
$$
B_\phi[D] = \Tr \phi(D^2),
$$
with $\phi$ being an ``arbitrary'' positive function of the Dirac
operator~$D$.

Chamseddine and Connes' work on the universal bosonic functional has
two main parts. In the first one, they argue that $B_\phi$ has the
following asymptotic development:
$$
B_\phi[D/\La] \sim \sum_{n=0}^\infty f_n\, \La^{4-2n}\, a_n(D^2)
\as \La \to \infty,
\eqno (1.1)
$$
where the $a_n$ are the coefficients of the heat kernel
expansion~\cite{\Gilkey} for $D^2$ and
$f_0 = \int_0^\infty x \phi(x) \,dx$,
$f_1 = \int_0^\infty \phi(x) \,dx$, $f_2 = \phi(0)$,
$f_3 = -\phi'(0)$, and so on. Then they proceed to compute the
development for the $K$-cycle currently~\cite{\ConnesNCGR, \Cordelia}
associated to the Standard Model, indeed obtaining all terms in the
bosonic part of the action for the Standard Model, plus gravity, plus
some new ones. Their approach gives prima facie relations between the
parameters of the Standard Model, in terms of the cutoff parameter
$\La$, falling rather wide of the empirical mark. In the second part
of their paper, they enterprise to improve the situation by use of the
renormalization group flow equations~\cite{\Persephone}. This need not
concern us here.

Formula~\eq{1.1} can be given a quick derivation, by assuming that
$\phi$ is a Laplace transform. This condition, however, will almost
never met in practice. In order to see that the asymptotic development
of $B_\phi$ cannot be taken for granted, let us consider, as Kastler
and coworkers have done~\cite{\CarminatiIKS, \IochumKS} the
characteristic functions $\phi_\La := \chi_{[0,\La]}$. This looks
harmless enough, giving nothing but $N_{D^2}(\La^2)$, the counting
number of eigenvalues of $D^2$ below the level $\La^2$. However, it
has been known for a long time ---see for
instance~\cite{\HormanderS}--- that there is no asymptotic development
for the counting functional beyond the first term. Therefore
equation~\eq{1.1}, as it stands, is not applicable to that situation.

One of our aims in this paper is to decrypt the meaning of ``arbitrary
functional''; a related one is to put on a firm footing the
development~\eq{1.1}. Our contribution turns around the Ces\`aro
behaviour of distributions, and its relation with asymptotic analysis.
Most results are new, or seem ignored in the literature; the paper is
written with a pedagogical bent.

The article is organized as follows. Section~2 is the backbone of the
paper; there the Ces\`aro behaviour of distributions and Ces\`aro
summability of evaluations are examined. The distributional theory of
asymptotic expansions~\cite{\EstradaKbook} is summarized. The latter
is brought to bear by finding the essential equivalence between the
Ces\`aro behaviour and the parametric behaviour of distributions at
infinity. Also we prove that a distribution satisfies the moment
asymptotic expansion iff it belongs to $\K'$, the dual of the space of
Grossmann--Loupias--Stein operator symbols~\cite{\GrossmannLS}. These
results are new, having been obtained very recently by one of
us~\cite{RE, \EstradaC}. We try to enliven this somewhat technical
section with pertinent examples.

Next we consider elliptic, positive pseudodifferential operators; let
$H$ be one of those; the functional calculus for $H$ can be based on
the \textbf{spectral density}, formally written as $\dl(\la - H)$.
This is arguably a more basic object than the heat kernel, and its
study is very rewarding. In Section~3, we show that $\dl(\la - H)$ is
an operator-valued distribution in $\K'$. With that in hand, one can
proceed to give a meaning to the universal bosonic action for a very
wide class of functionals. Following some old ideas by
Fulling~\cite{\FullingA}, insufficiently exploited up to now, we
emphasize that the Ces\`aro behaviour of the spectral density for
differential operators is local, i.e., independent of the boundary
conditions. This is practical for computational purposes, as it
sometimes allows to replace an operator in question by a more
convenient local model.

In Section~4, we reach the heart of the matter: let $d_H(x,y;\la)$
denote the distributional kernel of $\dl(\la - H)$; a formula for
$d_H$ is given and immediately applied to compute the coefficients of
its asymptotic expansion on the diagonal, in terms of the
noncommutative residues~\cite{\Wodzicki} of certain powers of~$H$. We
hope to have clarified in the paper that the identification of the
higher Wodzicki terms is essentially a ``finite-part'' calculation.
The spectral density is actually a less singular object for operators
with continuous spectra than for operators with discrete spectra, and
all of the above applies to operators associated to noncompact
manifolds: for that purpose, taking account of locality, we work with
densities of noncommutative residues throughout. We go on to extend
Connes' trace theorem~\cite{\ConnesA} to noncompact $K$-cycles. The
case of generalized Laplacians is then treated within our procedure.

In the light of the preceding, the last two sections of the paper are
concerned, respectively, with the counting number and the heat kernel
expansions. The counting functional $N_H(\la)$ is treated mainly by
way of example. Then we reexamine the status of arbitrary smoothing
asymptotic expansions, in particular the Laplace-type expansions like
the Chamseddine--Connes Ansatz. We point out conditions for the
expansions to be valid without qualification, and to be valid only in
the Ces\`aro sense. Also we exemplify circumstances under which the
formal Laplace-type expansion does not say anything about the true
asymptotic development. The Chamseddine--Connes expansion is derived
and reinterpreted.

\section 2. Ces\`aro computability of distributions

Besides the standard spaces of test functions and distributions, the
space $\K$ first introduced in~\cite{\GrossmannLS} and its dual $\K'$
play a central role in our considerations. Familiarity with the
properties of $\K$ and $\K'$ and with some of their elements will be
convenient. For all general matters in distribution theory, we refer
to~\cite{\GelfandShilov}.

As our interest is mainly in spectral theory, we consider
Grossmann--Loupias--Stein symbols in one variable, almost exclusively.
A smooth function $\phi$ of a real variable belongs to $\K_\ga$ for a
real constant $\ga$ if $\phi^{(k)}(x) = O(|x|^{\ga - k})$ as
$|x| \to \infty$, for each $k \in \N$. A topology for $\K_\ga$ is
generated by seminorms $\|\phi\|_{k,\ga}
 = \sup_{x\in\R} \set{\max(1, |x|^{k-\ga}) \,|\phi^{(k)}(x)|}$, and so
$\K_\ga \hookrightarrow \K_{\ga'}$ if $\ga \leq \ga'$. Notice that
$\phi^{(k)} \in \K_{\ga-k}$ if $\phi \in \K_\ga$. The space $\K$ is
the inductive limit of the spaces $\K_\ga$ as $\ga \to \infty$.

Since every polynomial is in~$\K$, a distribution $f \in \K'$ has
\textbf{moments}
$$
\mu_n := \<f(x), x^n>, \qquad n \in \N
$$
of all orders; this is an indication that $f$ decays rapidly at
infinity in some sense.

Denote by $\D'_0(\T)$ the space of periodic distributions with zero
mean. They constitute a first class of examples: if $f \in \D'_0(\T)$,
then, for $n$ suitably large, the periodic primitive with zero mean
$f_n$ of $f$ of order $n$ is continuous and defines the evaluation of
$f$ at $\phi \in \K$ by a convergent integral:
$$
\<f(x), \phi(x)> = (-1)^n \<f_n(x), \phi^{(n)}(x)>.
$$
Note that in this case all the moments are zero.

The algebra $\K$ is normal (i.e., $\SS$ is dense in $\K$) and is a
subalgebra of the multiplier algebras $\Oh_M$, $\M$ of $\SS$,
respectively for the ordinary product and the Moyal star
product~\cite{\Amalthea}. Other properties of $\K$ and $\K'$ will be
invoked opportunely. The usefulness of $\K$ in phase-space Quantum
Mechanics lies in the similitude of behaviour of the ordinary and the
Moyal product, when applied to elements of $\K$. The link between both
appearances of $\K$ is still mysterious to us.

\smallskip

The natural method of studying generalized functions at infinity is by
considering the parametric behaviour. The \textbf{moment asymptotic
expansion} of a distribution~\cite{\EstradaKbook} is given by
$$
f(\la x) \sim \sum_{k=0}^\infty
\frac{(-1)^k \mu_k\,\dl^{(k)}(x)}{k!\,\la^{k+1}} \as \la \to \infty.
\eqno (2.1)
$$
The interpretation of this formula is in the distributional sense, to
wit
$$
\<f(\la x), \phi(x)>
= \sum_{k=0}^N \frac{\mu_k \phi^{(k)}(0)}{k!\,\la^{k+1}}
+ O\Bigl(\frac{1}{\la^{N + 2}}\Bigr) \as \la \to \infty,
$$
for each $\phi$ in an appropriate space of test functions. Such an
expansion holds only for distributions that decay rapidly at infinity,
in a sense soon to be made completely precise; it certainly does not
hold for all tempered distributions, as their moments do not generally
exist. Distributions endowed with moment asymptotic expansions are
said to be ``distributionally small at infinity''. We are not happy
with this terminology and invite suggestions to improve it.

\smallskip

On the other hand, the classical analysis~\cite{\Hardy} notion of
Ces\`aro or Riesz means of series and integrals admits a
generalization to the theory of distributions, that we intend to
exploit in this paper. It turns out that Ces\`aro limits and
``distributional'' ones are essentially equivalent; this will enable
us to apply the simpler ideas of parametric analysis to complicated
averaging schemes.

We begin now in earnest by introducing the basic concept of Ces\`aro
behaviour of the distributions; justification will follow shortly.
Assume $f \in \D'(\R)$, $\b \in \R \setminus \{-1,-2,\dots\}$.

\declare Definition 2.1.
We say that $f$ is of order $x^\b$ at infinity, in the Ces\`aro sense,
and write
$$
f(x) = O(x^\b) \ala C \as x \to \infty,
$$
if there exists $N \in \N$, a primitive $f_N$ of $f$ of order $N$ and
a polynomial $p$ of degree at most $N - 1$, such that $f_N$ is locally
integrable for $x$~large and the relation
$$
f_N(x) = p(x) + O(x^{N + \b}) \as x \to \infty
\eqno (2.2)
$$
holds in the ordinary sense.

The relation $f(x) = o(x^\b) \ala C$ is defined similarly. The
notation $(C,N)$ can be used if one needs to be more specific; if an
order relation holds $(C,N)$ for some $N$, it also holds $(C,M)$ for
all $M > N$. The assumption $\b \neq -1, -2,\dots$ is provisionally
made in order to avoid dealing with the primitives of $x^{-1}$,
$x^{-2}$ and such (see Section~6 for the general case). If $\b > -1$,
the polynomial $p$ is arbitrary and thus irrelevant. We shall suppose
when needed that our distributions have bounded support, say, on the
left. In that case, we denote by $I[f]$ the first order primitive of
$f$ with support bounded on the left. When $f$ is locally integrable,
then,
$$
I[f](x) = \int_{-\infty}^x f(t)\, dt.
$$

The notation
$$
f(x) = o(x^{-\infty}) \ala C \as x \to \infty
$$
will mean $f(x) = O(x^\b) \ala C$ for every $\b$.

For the proof of the following workhorse proposition we refer
to~\cite{\EstradaC}.

\proclaim Lemma 2.1.
{\rm(a)}\enspace
Let $f\in\D'$ such that
$$
f(x) = O(x^\b) \ala{C,N} \as x \to \infty.
$$
Then for $k = 1, 2, 3,\dots$ we have:
$$
f^{(k)}(x) = O(x^{\b - k}) \ala{C,N+k} \as x \to \infty.
$$
{\rm(b)}\enspace
Let $f \in \D'$ such that
$$
f(x) = O(x^\b) \ala C \as x \to \infty,
$$
and let $\a \in \R$. Provided that $\a + \b$ is not a negative
integer, we have:
$$
x^\a f(x) = O(x^{\a + \b}) \ala C \as x \to \infty.
\eqno\qed
$$

\declare Definition 2.2.
We write $\lim_{x\to\infty} f(x) = L \ala C$ when
$f(x) = L + o(1) \ala C$ as $x \to \infty$. That is,
$\lim_{x\to\infty} f(x) = L \ala{C,k}$ when
$f_k(x)\,k!/x^k = L + o(1)$, for $f_k$ a primitive of order $k$
of~$f$.

For example, if $f$ is periodic with zero mean value, there exists
$n \in \N$ and a con\-tinuous (thus bounded) periodic function $f_n$
with zero mean value such that $f_n^{(n)} = f$; then clearly
$$
f(x) = o(x^{-\infty}) \ala C \as x \to \infty,
$$
a fact that yields, for $f$ periodic with mean value $a_0$:
$$
\lim_{x\to\infty} f(x) = a_0 \ala C.
$$

Let $f \in \D'$ be a distribution with support bounded on the left and
let $\phi$ be a smooth function. The following is a key concept of the
theory.

\declare Definition 2.3.
We say that the $\<f(x), \phi(x)>$ has the value $L$ in the Ces\`aro
sense, and write
$$
\<f(x), \phi(x)> = L \ala C
$$
if there is a primitive $I[g]$ for the distribution
$g(x) = f(x)\phi(x)$, satisfying
$$
\lim_{x\to\infty} I[g](x) = L \ala C \as x \to \infty.
$$

A similar definition applies when $f$ has support bounded on the
right. If $f$ is an arbitrary distribution, let $f = f_1 + f_2$ be a
decomposition of $f$, where $f_1$ has support bounded on the left and
$f_2$ has support bounded on the right. Then we say that
$\<f(x), \phi(x)> = L\ala C$ if both $\<f_i(x), \phi(x)> = L_i \ala C$
exist for $i = 1,2$ and $L = L_1 + L_2$: this definition is seen to be
independent of the decomposition.

For instance, let $f$ be a periodic distribution of zero mean and let
$f_1,f_2,\dots,f_{n+1}$ denote the periodic primitives with zero mean
of $f$, up to the order $n + 1$. Then
$$
x^n f_1(x) - nx^{n-1} f_2(x) + n(n-1)x^{n-2} f_3(x) - \cdots
+ (-1)^n n!\, f_{n+1}(x)
$$
is a first order primitive of $x^n f(x)$, and since
$f_i(x) = o(x^{-\infty})\ala C$ for $i = 1,\dots,n$ as $x \to \infty$,
it follows that
$$
\<f(x), x^k> = 0 \ala C
$$
for all $k \in \N$.

\smallskip

To perceive the point of our hitherto abstract definitions, it is
worthwhile to recall here briefly the classical theory~\cite{\Hardy}.
Let $\{a_n\}_{n=1}^\infty$ be a sequence of real or complex numbers.
Often it has no limit, but the sequence of averages
$H^{(1)}_n := (a_1 +\cdots+ a_n)/n$ does. Then people write
$$
\lim_{n\to\infty}a_n = L \ala{C,1}.
$$

If $H^{(1)}_n$ still does not have a limit, then one may apply the
averaging procedure again and again, hoping that eventually a limit
will be obtained. There are two main procedures to perform such higher
order averages: the H\"older means and the Ces\`aro means. The
H\"older means are single-mindedly constructed as
$$
H^{(k)}_n := \frac{H^{(k-1)}_1 + \cdots + H^{(k-1)}_n}{n}
$$
and $\lim_{n\to\infty} H^{(k)}_n = L$ is written
$$
\lim_{n\to\infty} a_n = L \ala{H,k}.
$$

The properly named Ces\`aro means are defined as follows: let
$A_n^{(0)} := a_n$ and define recursively
$A_n^{(k)} = A_1^{(k-1)} +\cdots+ A_n^{(k-1)}$. If
$\lim_{n\to\infty} k!\,A_n^{(k)}/n^k = L$, we write
$$
\lim_{n\to\infty}a_n = L \ala{C,k},
$$
so that the $(C,1)$ and the $(H,1)$ notions are identical. The
Ces\`aro limits have nicer ana\-lytical properties. The good news, at
any rate, is that both procedures are equivalent:
$$
\lim_{n\to\infty} a_n = L \ala{C,k}
\iff \lim_{n\to\infty} a_n = L \ala{H,k}.
$$
One uses the simpler notation $\lim_{n\to\infty} a_n = L \ala C$ if
$\lim_{n\to\infty} a_n = L \ala{C,k}$ for some $k\in\N$.

A third averaging procedure is equivalent to Ces\`aro's, the so-called
Riesz typical means. For real $\mu$, one writes
$$
\lim_{n\to\infty} a_n = L \ala{R,k,n}
$$
if
$$
\lim_{\mu\to\infty} \frac{1}{\mu}
\sum_{n\leq\mu} \Bigl(1 - \frac{n}{\mu}\Bigr)^{k-1} a_n = L.
$$
Riesz originally studied this formula for integral $\mu$, but the
means have more desirable properties with $\mu$ real. Now, one may
study the summability of a series $\sum_{n=1}^\infty a_n$ by studying
the generalized function of a real variable
$f(x) = \sum_{n=1}^\infty a_n \,\dl(x - n)$. The definition of
Ces\`aro limits of distributions is tailored in such a way that
$\<f,1> \ala C$ and $\sum_{n=1}^\infty a_n \ala C$ coincide: a
primitive of order $k$ of $\sum_{n=1}^\infty a_n \,\dl(x - n)$ is
given by $f_k(x) = \sum_{n\leq x} (x - n)^{k-1} a_n/(k - 1)!$ Note
that one could consider distributions of the form
$h(x) = \sum_{n=1}^\infty a_n \,\dl(x - p_n)$, with
$p_n \uparrow \infty$; this gives rise to the $(R,k,p_n)$ means.

In summary, we have demonstrated the following equivalence.

\proclaim Theorem 2.2.
The evaluation
$$
\biggl< \sum_{n=1}^\infty a_n\,\dl(x - n), \phi(x) \biggr> = L \ala C
$$
holds iff $\sum_{n=1}^\infty a_n \,\phi(n) = L$ in the Ces\`aro sense
of the theory of summability of series.
\qed

In the same vein:

\proclaim Theorem 2.3.
If $f$ is locally integrable and supported in $(a,\infty)$, then
$$
\<f(x), \phi(x)> = L \ala C
$$
if and only if
$$
\int_a^\infty f(x)\phi(x) \,dx = L
$$
in the Ces\`aro sense of the theory of summability of integrals.
\qed

As shown below, if $f\in\K'$ and $\phi\in\K$, then the evaluation
$\<f(x), \phi(x)>$ is always $(C)$-sum\-ma\-ble. We pause an instant
to show by the example just how useful is the concept of Ces\`aro
computability of evaluations. An interesting periodic distribution is
the Dirac comb $\sum_{n=-\infty}^\infty \dl(x - n)$. Its mean value
is~$1$; therefore
$$
\sum_{n=-\infty}^\infty \dl(x - n) = 1 + f(x),
\eqno (2.3)
$$
with $f \in \D'_0(\T)$. The distributions
$$
\sum_{n=1}^\infty \dl(x - n) - H(x - 1), \qquad
\sum_{n=1}^\infty \dl(x - n) - H(x),
$$
where $H$ is the Heaviside function, belong to $\K'$. In effect, take
a function $\phi_1 \in \K$ such that $\phi_1(x) = 1$ for $x > 1/2$,
$\phi_1(x) = 0$ for $x < 1/4$. Then
$\phi_1(x) \bigl( \sum_{n=-\infty}^\infty \dl(x - n) - 1 \bigr)$ only
differs from $\sum_{n=1}^\infty \dl(x - n) - H(x - 1)$ or
$\sum_{n=1}^\infty \dl(x - n) - H(x)$ by a distribution of compact
support.

It follows that the evaluation
$$
\biggl< \sum_{n=1}^\infty \dl(x - n) - H(x - 1), \phi(x) \biggr>
= \sum_{n=1}^\infty \phi(n) - \int_1^\infty \phi(x) \,dx
$$
is Ces\`aro summable whenever $\phi\in\K$. Now, $x^\a$ does not belong
to $\K$ unless $\a \in \N$, but the previous argument, using
$\phi_\a(x) = \phi_1(x)\,x^\a$, allows us to conclude that the
evaluation
$$
Z(\a) := \biggl< \sum_{n=1}^\infty \dl(x - n) - H(x - 1), x^\a \biggr>
$$
is $(C)$-summable for any $\a \in \C$. Also, $Z(\a)$ is an entire
function of $\a$, since $\phi_\a$ is. We find a formula for $Z(\a)$ by
observing that if $\Re\,\a < -1$ then the evaluation is given by the
difference of a series and an integral, so that
$$
Z(\a) = \sum_{n=1}^\infty n^\a - \int_1^\infty x^\a
= \zeta(-\a) + \frac{1}{\a + 1}, \quad \Re \a < -1.
$$

We have learned a simple proof that Riemann's zeta function is
analytic in $\C \setminus \{1\}$, with residue at $s = 1$ equal
to~$1$, and one realizes that the evaluation of the $\zeta$ function
can be done by Ces\`aro means (it is only because the zeta function is
the outcome of a regularization process that it is useful for
renormalization in quantum field theory). The evaluation
$\bigl< \sum_{n=1}^\infty \dl(x - n) - H(x), x^\a \bigr>$ is slightly
more involved. However, we may write~\cite{\EstradaKfp}:
$$
\biggl< \sum_{n=1}^\infty \dl(x - n) - H(x), x^\a \biggr>
:= Z(\a) - \Fp\int_0^1 x^\a \,dx,
$$
where $\Fp$ stands for the Hadamard finite part of the integral. Now,
$$
\Fp\int_0^1 x^\a dx = \frac{1}{\a + 1}, \qquad \a \neq -1,
$$
therefore, if $\a \neq -1$,
$$
\zeta(-\a)
= \sum_{n=1}^\infty n^\a - \Fp\int_0^\infty x^\a \,dx \ala C,
$$
in the sense that
$$
\zeta(-\a) = \lim_{x\to\infty} \biggl( \sum_{n=1}^{\piso{x}} n^\a
- \Fp\int_0^x t^\a \,dt \biggr) \ala C.
$$
This formula gives a nice representation for $\zeta(\a)$ when
$\Re\,\a < 1$. For instance, $\zeta(0) = -1/2$ simply because the
fractional part $\fpart{x} = x - \piso{x}$ of $x$ is periodic of mean
$1/2$. For $\a = -1$:
$$
\zeta(-1) = \lim_{x\to\infty} \biggl(
\sum_{n=1}^{\piso{x}} n - \int_0^x t \,dt \biggr) = \lim_{x\to\infty}
\bigl( \thalf \piso{x}(\piso{x} + 1) - \thalf x^2 \bigr) \ala C;
$$
we find that
$$
\frac{(x - \fpart{x})(x - \fpart{x} + 1)}{2} - \frac{x^2}{2}
= \frac{\fpart{x}^2 - \fpart{x}}{2} + \frac{x(1 - 2\fpart{x})}{2}
= -\frac{1}{12} + o(x^{-\infty}) \ala C
$$
since $(1 - 2\fpart{x})$ and $(\fpart{x}^2 - \fpart{x} + 1/6)$ are
periodic of mean zero; we get $\zeta(-1) = -1/12$.

Also, the logarithm of the ``functional determinant'' can be obtained
by this method:
$$
\zeta'(0) = - \lim_{x\to\infty} \biggl(
\sum_{n=2}^{\piso{x}} \log n - \int_0^x \log t \,dt \biggr) \ala C,
$$
on using Lemma~2.1. Stirling's formula gives
$$
\eqalign{
x \log x - x &- \log(\piso{x}!)
= x \log x - x - (\piso{x} + \thalf) \log\piso{x} + \piso{x}
- \log \sqrt{2\pi} + O(x^{-1})
\cr
&= - x \log\Bigl(1 - \frac{\fpart{x}}{x}\Bigr) - \fpart{x}
+ (\fpart{x} - \thalf) \log\piso{x} - \thalf \log(2\pi) + O(x^{-1})
\cr
&= - \thalf \log(2\pi) + O(x^{-1}) \ala C
\cr}
$$
since $x \log(1 - x^{-1}\fpart{x}) + \fpart{x} = O(x^{-2})$ and
$(\fpart{x} - \thalf)$ is periodic of mean zero. {}From this it
follows that $\zeta'(0) = -\thalf \log(2\pi)$. This business of
Riemann's zeta function is not merely amusing; it will be useful
later.

\smallskip

We make ready for the main equivalence result.

\proclaim Theorem 2.4.
Let $f \in \D'$. If $\a > -1$ then
$$
f(x) = O(|x|^\a) \ala C \as x \to \pm\infty
\eqno (2.4)
$$
if and only if
$$
f(\la x) = O(\la^\a) \as \la \to \infty
\eqno (2.5)
$$
in the topology of $\D'$. If $-j-1 > \a > - j -2$ for some $j \in \N$,
then~\eq{2.4} holds if and only if there are constants
$\mu_0,\dots,\mu_j$ such that
$$
f(\la x) = \sum_{k=0}^j
\frac{(-1)^k \mu_k\,\dl^{(k)}(x)}{k!\,\la^{k+1}} + O(\la^\a)
$$
in the topology of $\D'$ as $\la \to \infty$.

\Proof.
We prove the theorem in the case $f$ has support bounded on the left.
The general case follows by using a decomposition $f = f_1 + f_2$,
where $f_1$ has support bounded on the left and $f_2$ has support
bounded on the right. First we have to clarify the meaning
of~\eq{2.5}. It is a weak or distributional relation: we write
$f(x,\la) = O(\la^\a)$ as $\la \to \infty$ whenever
$$
\<f(x,\la), \phi(x)> = O(\la^\a) \as \la \to \infty,
$$
for all $\phi \in \D$. Note that this yields
$$
\biggl< \pd{f(x,\la)}{x}, \phi(x) \biggr>
= - \<f(x,\la), \phi'(x)> = O(\la^\a).
$$
Now, if~\eq{2.5} holds, there exists $N$ such that the primitive of
order $N$ of $f(\la x)$, with respect to $x$, exists and is bounded by
$M\la^\a$, say for $|x|\leq 1$ and $\la \geq \la_0$. We have then a
primitive $f_N$ of order $N$ of $f(x)$, such that
$$
|f_N(\la x)| \leq M\la^{\a + N}, \qquad |x|\leq 1,\ \la \geq \la_0.
$$
Taking $x = 1$ and replacing $\la$ by $x$ we obtain
$$
|f_N(x)| \leq Mx^{\a + N}, \qquad x \geq \la_0,
$$
and thus
$$
f(x) = O(x^\a) \ala{C,N}, \as x \to \infty.
$$
Reciprocally, assume $\a > -1$ and $f(x) = O(x^\a) \ala{C,N}$, as
$x \to \infty$. Then, if $f_N$ is the (locally integrable for $x$
large) primitive of order $N$ of $f$ with support bounded on the left,
an obvious estimate gives $f_N(\la x) = O(\la^{\a + N})$, as
$\la\to\infty$, and on differentiating $N$ times with respect to $x$
one obtains $\la^N f(\la x) = O(\la^{\a + N})$, so that~\eq{2.5}
follows.

The case when $\a$ is nonintegral and less than $-1$ is more involved,
as one has to deal with the polynomial $p$ in~\eq{2.2}. Then one shows
that the moments
$$
\<f(x), x^k> = \mu_k \ala C
$$
up to a certain order exist, those being essentially the coefficients
of~$p$. For the gory details, we refer once again to~\cite{\EstradaC}.
\endproof

A characterization of the distributions that have a moment asymptotic
expansion follows.

\proclaim Theorem 2.5.
Let $f \in \D'$. Then the following are equivalent:
\item{\rm(a)} $f \in \K'$.
\item{\rm(b)} $f$ satisfies
$$
f(x) = o(|x|^{-\infty}) \ala C \as x \to \pm\infty.
$$
\item{\rm(c)} There exist constants $\mu_0,\mu_1,\mu_2,\dots$ such
that
$$
f(\la x) \sim \frac{\mu_0\,\dl(x)}{\la} - \frac{\mu_1\,\dl'(x)}{\la^2}
+ \frac{\mu_2\,\dl''(x)}{2!\,\la^3} - \cdots \as \la \to \infty
$$
in the weak sense.

\Proof.
It is proven in~\cite{\EstradaKbook} that the elements of $\K'$
satisfy the moment asymptotic expansion. For the converse, it is
enough, as customary, to consider distributions with support bounded
on one side. We show that if (b) holds, then $f\in K'_\ga$ for all
$\ga$. {}From the hypothesis it follows that
$f(x) = O(x^{-\ga-2}) \ala C$ as $x\to\infty$. Thus, for a certain
$n$, the $n$-th order primitive $f_n$ of~$f$ with support bounded on
one side is locally integrable and satisfies
$f_n(x) = p(x) + O(x^{-\ga-2+n})$ as $x\to\infty$, where the
polynomial $p$ has degree at most $n - 1$. We conjure up a compactly
supported continuous function $g$ whose moments of order up to $n - 1$
coincide with those of~$f$. If $g_n$ is the primitive of order $n$ of
$g$ with support bounded on the left, then
$f_n(x) - g_n(x) = O(x^{-\ga-2+n})$. If $\phi \in \K_{\ga-n}$, the
integral $\int_{-\infty}^\infty (f_n(x) - g_n(x)) \phi(x)\,dx$
converges. Hence $f = (f_n - g_n)^{(n)} + g \in \K'_\ga$. The rest is
clear.
\endproof

We get at once a powerful computational method for duality
evaluations.

\proclaim Corollary 2.6.
If $f \in \K'$ and $\phi \in \K$, the evaluation $\<f(x), \phi(x)>$
is Ces\`aro summable.

\Proof.
It is enough to check for $\phi = 1$. But, according to the previous
Theorem, if $f\in\K'$ then $f(x) = o(x^{-\infty}) \ala C$ as
$x \to \infty$. By the proof of Theorem~2.4, $\<f(x), 1>$ is
$(C)$-summable.
\endproof

Fourier transforms are defined by duality and, in general, if
$f \in \SS'$, we cannot make sense of $\hat f(u)$ because the
evaluation $\<e^{ixu}, f(x)>$ is not defined. However, if
$\phi \in \K$ and $u \neq 0$, Corollary~2.6 guarantees that the
Ces\`aro-sense evaluation $\<e^{ixu}, \phi(x)> \ala C$ is well
defined. Thus
$$
\hat\phi(u) = \<e^{ixu}, \phi(x)> \ala C
\sepword{when} \phi \in \K, \ u \neq 0.
$$
It is clear that $\widehat\K \subset \K'$; this follows also from
Proposition 4 of~\cite{\GrossmannLS}.

Note as well that the moments of $f \in \K'$ are $(C)$-summable. The
converse is true:

\proclaim Theorem 2.7.
Let $f \in \D'$. If all the moments $\<f(x), x^n> =\mu_n \ala C$ exist
for $n \in \N$, then $f \in \K'$.

For the easy proof, we refer to~\cite{\EstradaC}.
\endproof

It is clearly worthwhile to characterize spaces of distributions in
terms of their Ces\`aro behaviour. Particularly important is the
characterization of tempered distributions:

\proclaim Theorem 2.8.
Let $f \in \D'$. Then the following statements are equivalent:
\item{\rm(a)} $f$ is a tempered distribution.
\item{\rm(b)} There exists $\a \in \R$ such that
$$
f(\la x) = O(\la^\a), \as \la \to \infty
$$
in the weak sense.
\item{\rm(c)} There exists $\a \in \R$ and $k \in \N$ such that
$$
f^{(k)}(x) = O(|x|^{\a - k}) \ala C \as x \to \infty.
$$

\Proof.
Again, it is enough to consider the case when $f$ has support bounded
on one side. It is well known that if $f \in \SS'$ then there is a
primitive $F$ of some order $N$ of slow growth at infinity; it follows
that $f(x) = O(|x|^\a) \ala C$. The rest is clear, in view of the
equivalence theorem~2.4 and the fact that distributional order
relations can be differentiated at will.
\endproof

We finish by giving several estimates that we will need later. The
first one is just a rewording of the properties of the
distribution~\eq{2.3}.

\proclaim Lemma 2.9.
If $g \in \K$ and if $\int_{-\infty}^\infty g(x)\,dx$ is defined, then
$$
\sum_{n=-\infty}^\infty g(n\eps)
= \frac{1}{\eps} \int_{-\infty}^\infty g(x)\,dx + o(\eps^\infty)
\as \eps \downarrow 0.
\eqno\qed
$$

By the same token:

\proclaim Lemma 2.10.
If $g \in \K(\R^n)$ and if $\int_{\R^n} g(x)\,dx$ is defined, then
$$
\sum_{k\in\Z^n} g(k\eps)
= \eps^{-n} \int_{\R^n} g(x)\,dx + o(\eps^\infty)
\as \eps \downarrow 0.
\eqno\qed
$$

\proclaim Lemma 2.11.
If $g \in \K$ and if $\int_0^\infty g(x)\,dx$ is defined, then
$$
\sum_{n=1}^\infty g(n\eps) = \frac{1}{\eps} \int_0^\infty g(x)\,dx
+ \sum_{n=0}^\infty \frac{\zeta(-n) g^{(n)}(0)}{n!} \eps^n +
o(\eps^\infty) \as \eps \downarrow 0.
$$

\Proof.
This follows from the zeta function example.
\endproof

(Results of this type were used to prove some formulas by Ramanujan
in~\cite{\EstradaKbook}.)

\section 3. Spectral densities

Let $\H$ be a concrete Hilbert space, the space of square integrable
sections of an Euclidean vector bundle over a Riemannian manifold $M$,
and let $H$ be an elliptic positive selfadjoint pseudodifferential
operator on $\H$, with domain $\X$. We consider the derivative, in the
distributional sense, of the spectral family of projectors $E_H(\la)$
associated to $H$:
$$
d_H(\la) := \frac{dE_H(\la)}{d\la}.
$$
For instance, if $H$ is defined on a compact manifold, and
$0 < \la_1 \leq \la_2 \leq \cdots$ is the complete set of its
eigenvalues, with orthonormal basis of eigenfunctions $u_j$, the
kernel of the spectral family is given by~\cite{\HormanderIII}:
$$
E_H(\la) := \sum_{\la_j\leq\la} \ketbra{u_j}{u_j},
$$
and the derivative is
$$
d_H(\la) := \sum_j \ketbra{u_j}{u_j} \,\dl(\la - \la_j).
$$

This \textbf{spectral density} is a distribution with values in
$\L(\X,\H)$. The defining properties of $E(\la)$:
$$
I = \int_{-\infty}^\infty dE(\la), \qquad
H = \int_{-\infty}^\infty \la \,dE(\la)
$$
(in the weak sense) become, in the language of the previous section:
$$
I = \<d_H(\la), 1>, \qquad H = \<d_H(\la), \la>.
$$

The spectral density is used to construct the functional calculus
for~$H$. Indeed, we can define $\phi(H)$ whenever $f$ is a
distribution such that the evaluation $\<d_H(\la), f(\la)>$ makes
sense, by
$$
\phi(H) := \<d_H(\la), \phi(\la)>,
$$
with domain the subspace of the $x \in \H$ for which the evaluation
$\bigl< \braket{y}{d_H(\la)x}, \phi(\la) \bigr>_\la$ is defined for
all $y \in \H$.

Especially, one is able to deal with the ``zeta operator":
$$
H^{-s} := \<d_H(\la), \la^{-s}>,
\eqno (3.1)
$$
(for $0 \notin \spec H$), the heat operator:
$$
e^{-tH} := \<d_H(\la), e^{-t\la}>,  \quad  t > 0
\eqno (3.2)
$$
and the unitary group of $H$, which is just the Fourier transform of
the spectral density:
$$
U_H(t) := \<d_H(\la), e^{-it\la}>.
\eqno (3.3)
$$
The useful symbolic formula
$$
d_H(\la) = \dl(\la - H)
$$
recommends itself, and we shall employ it from now on.

We want to study the asymptotic behaviour of $\dl(\la - H)$. Let
$\X_n$ be the domain of $H^n$ and let
$\X_\infty := \bigcap_{n=1}^\infty \X_n$. The fact that $\X_\infty$ is
dense has, in view of the theory of Section~2, momentous consequences.
We have
$$
H^n = \<\dl(\la - H), \la^n>
$$
in the space $\L(\X_\infty,\H)$. Hence, $\dl(\la - H)$ belongs to the
space $\K'(\R, \L(\X_\infty,\H))$. Therefore the moment asymptotic
expansion holds:
$$
\dl(\la\sigma - H) \sim \sum_{n=0}^\infty
\frac{(-1)^n H^n \,\dl^{(n)}(\la)}{n!\,\sigma^{n+1}}
\as \sigma \to \infty,
$$
and $\dl(\la - H)$ vanishes to infinite order at infinity in the
Ces\`aro sense:
$$
\dl(\la - H) = o(|\la|^{-\infty}) \ala C \as |\la| \to \infty.
$$
Of course, the last formula is trivial when $H$ is bounded.

The space $\D(M)$ of test functions is a subspace of $\X_\infty$. We
can then realize the spectral density by an associated kernel
$d_H(x,y;\la)$, an element of $\D'(\R, \D'(M \x M))$. Ellipticity
actually implies that $d_H(x,y;\la)$ is smooth in $(x,y)$. The
expansion
$$
d_H(x,y; \la\sigma) \sim \sum_{n=0}^\infty \frac{(-1)^n
(H^n\dl)(x-y)\,\dl^{(n)}(\la\sigma)}{n!\,\sigma^{n+1}}
\as \sigma \to \infty
$$
holds in principle in the space $\D'(\R, \D'(M \x M))$. We also get
$$
d_H(x,y; \la) = o(|\la|^{-\infty}) \ala C \as |\la| \to \infty
\eqno (3.4)
$$
in the space $\D'(M \x M)$. Equation~\eq{3.4} is the mother of all
incoherence principles. For instance, passing to the primitive with
respect to $\la$, for an elliptic operator on a compact manifold with
eigenfunctions $\psi_n$, $n\in \N$, one concludes:
$$
\sum_{\la_n \leq \la}\bar\psi_n(x)\psi_n(y)
= o(|\la|^{-\infty}) \ala C \as |\la| \to \infty,
$$
for $x \neq y$, which is Carleman's incoherence
relation~\cite{\Carleman}.

It should be clear that the expansions cannot hold pointwise in both
variables $x$ and $y$, since we cannot set $x = y$ in the distribution
$\dl(x - y)$. In fact, our interest in this paper lies in the
coincidence limit $d_H(x,x;\la)$, which is \textbf{not}
distributionally small. However, it is proven in~\cite{\EstradaF}
that, away from the diagonal of $M \x M$, the expansions are valid in
the sense of uniform convergence of all derivatives on compacta. On
the other hand, if $H_1$ and $H_2$ are two pseudodifferential
operators whose difference over an open subset $U$ of $M$ is a
smoothing operator, and if $d_1(x,y;\la)$ and $d_2(x,y;\la)$ are the
corresponding spectral densities, then~\cite{\EstradaF}:
$$
d_1(x,y; \sigma\la) = d_2(x,y; \sigma\la) + o(\sigma^{-\infty})
\as \sigma \to \infty
$$
in $\D'(U \x U)$. Also, it can be shown that
$$
d_1(x,y;\la) = d_2(x,y;\la) + o(\la^{-\infty}) \ala C
\as \la \to \infty
$$
uniformly on compacts of $U \x U$, even at the diagonal.

We exemplify the reported behaviour with the simplest possible
examples. Let $H$ denote first the Laplacian on the real line. Its
spectral density is
$$
d_H(x,y;\la)
= \frac{1}{2\pi\sqrt{\la}} \cos\bigl(\sqrt{\la}(x - y)\bigr)
$$
and therefore it is clear that $d_H(x,x;\la)$ is not distributionally
small, but rather
$$
d_H(x,x;\la) = \frac{1}{2\pi\sqrt{\la}} + o(\la^{-\infty}) \ala C
\as \la \to \infty.
$$

Let $H$ denote now the Laplacian on the circle; the eigenvalues are
$\la_n = n^2$, $n = 0,1,2,\dots$, with multiplicity~$2$ from $n = 1$
on, with normalized eigenfunctions
$\psi_n^\pm(x) = (2\pi)^{-1/2}\, e^{\pm inx}$. Therefore
$$
d_H(x,y;\la) = \frac{1}{2\pi} \Bigl( \dl(\la)
+ 2 \sum_{n=1}^\infty \cos n(x - y) \,\dl(\la - n^2) \Bigr).
$$
Then
$$
\frac{1}{2\pi} \Bigl( \dl(\la\sigma)
+ 2\sum_{n=1}^\infty \cos n(x - y) \,\dl(\la\sigma - n^2) \Bigr)
\sim \sum_{j=0}^\infty
\frac{\dl^{(2j)}(x - y) \,\dl^{(j)}(\la)}{j!\,\sigma^{j+1}}
\as \sigma \to \infty
$$
in $\D'(\R,\D'(\Sf^1 \x \Sf^1))$, while
$$
\frac{1}{2\pi} \Bigl( \dl(\la)
+ 2\sum_{n=1}^\infty \cos n(x - y)\,\dl(\la - n^2) \Bigr)
= o(\la^{-\infty}) \ala C \as \la \to \infty
$$
if $x$ and $y$ are fixed, $x \neq y$.

On the other hand,
$$
d_H(x,x;\la) = \frac{1}{2\pi}
\Bigl( \dl(\la) + 2\sum_{n=1}^\infty \,\dl(\la - n^2) \Bigr)
$$
does not belong to $\K'(\R, C^\infty(M))$. For the first time in this
paper, but not the last, we have to find out what the Ces\`aro
behaviour of a given spectral kernel is. We shall have recourse to a
variety of tricks. For now, applying Lemma~2.11 to
$g(x) := \phi(x^2)$, for $\phi$ a Schwartz function, say, we get:
$$
\sum_{n=1}^\infty \phi(\eps n^2)
= \frac{1}{2\sqrt{\eps}} \int_0^\infty x^{-1/2} \phi(x)\,dx
- \thalf\phi(0) + o(\eps^\infty) \as \eps \downarrow 0.
$$
It is then clear that
$$
d_H(x,x;\la) = \frac{1}{2\pi\sqrt{\la}} + o(\la^{-\infty}) \ala C
\as \la \to \infty.
$$
and it is also immediately clear that the distributional and Ces\`aro
behaviour of the spectral density and its kernel are \textit{exactly
the same} as in the previous example. That the manifold be compact or
not and the spectrum be discrete or continuous is immaterial for that
purpose. If we seek a boundary problem for the Laplacian, say on a
bounded interval of the line, we obtain still the same kind of
behaviour (off the boundary, where a sharp change takes place). Note
also the estimate:
$$
\sum_{\pm;\,\la_n\leq\la}|\psi_n^\pm(x)|^2
\sim \frac{\sqrt{\la}}{\pi} \ala C \as \la \to \infty.
$$

As an aside, we turn before closing this section to the functional
calculus formulas and compare~\eq{3.2} with~\eq{3.3}. Obviously
$e^{-t(\cdot)}$ has an extension belonging to $\K$, so there is no
difficulty in giving a meaning to the heat operator. Also, as we shall
see in Section~6, it is comparatively easy to study the asymptotic
development of the corresponding Green function as $t \downarrow 0$.
One of the motivations of the present approach to spectral asymptotics
is to define a sense for expansions of Schr\"odinger propagators and
the like, that do not possess a ``true'' asymptotic expansion.

Such an approach can be based in the following idea: Theorem~2.8
points to a rough duality between $\K'$ and $\SS'$. Let
$g \in \SS'(\R)$ and find $\a$ so that $g(\la x) = O(\la^\a)$ weakly
as $\la \to \infty$. For any $\phi \in \SS(\R)$, the function $\Phi$
defined by
$$
\Phi(x) := \<g(\la x), \phi(\la)>_\la
$$
is smooth for $x \neq 0$ since
$\Phi(x) = |x|^{-1} \<g(\la), \phi(\la x^{-1})>_\la$, and satisfies
$$
\Phi^{(n)}(x) = O(|x|^{\a-n}) \as |x| \to \infty.
$$
Therefore, if $f \in \K'$ with $0 \notin \supp f$, we can define
$\bigl< f(x), g(\la x) \bigr>_x$ as a tempered distribution.

When $0 \in \supp f$, we need to ascertain independently smoothness of
$\Phi$ at the origin. It turns out that, for this purpose, it is
enough to demand distributional smoothness of $g$, i.e., the existence
of the \textit{distributional} values $g^{(n)}(0)$, in the sense
of~\cite{\Lojasiewicz}, for $n = 0,1,2,\dots$. Then $g(tH)$ admits a
\textit{distributional} expansion in $\L(\X_\infty,\H)$ as
$t \downarrow 0$. This can eventually lead to a proper treatment of
some questions in quantum field theory. We say no more here and refer
instead to the forthcoming~\cite{\EstradaF}. In Section 6 of this
paper, results will be stated for $g$ belonging to $\SS(\R)$; for the
rest of the paper we will venture outside safe territory only in
examples.

\section 4. The Ces\`aro asymptotic development of $d_H(x,x;\lambda)$

In this section we obtain the asymptotic expansion for the coincidence
limits of spectral density kernels. We are fortified with the results
of the previous section, implying that the Ces\`aro behaviour of the
spectral density of pseudodifferential operators is a local matter.

Let $A$ be any pseudodifferential operator of order a positive integer
$d$, with complete symbol $\sigma(A)$, on the Riemannian manifold $M$.
To simplify the discussion, we consider only operators acting on
scalars; the treatment of matrix-valued symbols presents no further
difficulty. The noncommutative or Wodzicki \textbf{residue} of $A$ is
defined by integrating (the trace of) the partial symbol
$\sigma_{-n}(A)(x,\xi)$ of order~$-n$ over the cosphere bundle
$\set{(x,\xi) : |\xi| = 1}$:
$$
\Wres A := \int_M \int_{\Sf^{n-1}} \sigma_{-n}(A)(x,\xi) \,d\xi \,dx.
$$
Here $dx$ denotes the canonical volume element on $M$. If $M$ is not
compact, $\Wres A$ may not exist, but there always exists the local
density of the residue $\int_{\Sf^{n-1}}\sigma_{-n}(A)(x,\om) \,d\om$,
that we denote by $\wres A(x)$.

We recall that
$$
\sigma(AB) - \sigma(A)\sigma(B) \sim \sum_{|\a|>0}
\frac{(-i)^{|\a|}}{\a!} \del_\xi^\a \sigma(A) \del_x^\a\sigma(B).
$$

The kernel $k_A$ of $A$ is by definition:
$$
k_A(x,y)
:= (2\pi)^{-n} \bigl< e^{i(x-y)\.\xi}, \sigma(A)(x,\xi) \bigr>_\xi.
$$
In particular, on the diagonal:
$$
k_A(x,x) := (2\pi)^{-n} \bigl< 1, \sigma(A)(x,\xi) \bigr>_\xi.
\eqno (4.1)
$$

In order to figure out the symbol for a spectral density, we start by
considering (the selfadjoint extension of) an elliptic operator $H$
with constant coefficients. In this case $\sigma(H^n) = \sigma(H)^n$
and we assert:
$$
\sigma\bigl( \dl(\la - H) \bigr) = \dl(\la - \sigma(H)),
$$
justified by the identities:
$$
\int \la^n \dl(\la - \sigma(H)) \,d\la = \sigma(H^n), \qquad
\la = 0, 1, 2,\dots
$$

In the general case of nonconstant coefficients, we make the Ansatz
that:
$$
\sigma\bigl( \dl(\la - H) \bigr)
\sim \dl(\la - \sigma(H)) - q_1\,\dl'(\la - \sigma(H))
+ q_2\,\dl''(\la - \sigma(H)) - q_3\,\dl'''(\la - \sigma(H)) + \cdots
\eqno (4.2)
$$
in the Ces\`aro sense. Computation of
$\int \la^n \sigma(\dl(\la - H))\,d\la$ for $\la = 0,1,2,\dots$ then
gives
$$
q_1 = 0; \quad
q_2 = \thalf\bigl(\sigma(H^2) - \sigma(H)^2\bigr); \quad
q_3 = \tfrac{1}{6}\bigl( \sigma(H^3) - 3\sigma(H^2)\sigma(H)
+ 2\sigma(H)^3 \bigr),
\eqno (4.3)
$$
and so on. This development, it turns out, gives ever lower powers of
$\la$ in the asymptotic expansion of $\sigma(\dl(\la - H))$.

We are interested in explicit formulas for the Ces\'aro asymptotic
development of the coincidence limit for the kernel of a positive
operator $H$ as $\la \to \infty$. {}From~\eq{4.1} and~\eq{4.2} with
$p := \sigma(H)$, we get
$$
d_H(x,x;\la) \sim (2\pi)^{-n} \bigl<1, \dl(\la - p(x,\xi))
+ q_2(x,\xi)\,\dl''(\la - p(x,\xi)) - \cdots \bigr>_\xi \ala C.
$$
In polar coordinates on the cotangent fibres, $\xi = |\xi|\om$ with
$|\om| = 1$, this becomes
$$
(2\pi)^{-n} \int_{|\om|=1} \! d\om\, \bigl< |\xi|^{n-1},
\dl(\la - p(x, |\xi|\om))
+ q_2(x, |\xi|\om)\,\dl''(\la - p(x,|\xi|\om))
- \cdots \bigr>_{|\xi|}.
$$
Hence, if we denote by $|\xi|(x,\om;\la)$ the positive solution of the
equation $p(x,|\xi|\om) = \la$, we need to compute:
$$
(2\pi)^{-n} \int_{\Sf^{n-1}} \! d\om\,
\frac{|\xi|^{n-1}(x,\om;\la) + \frac{\del^2}{\del\la^2}
\bigl( q_2(x,|\xi|(x,\om;\la)\om)|\xi|^{n-1}(x,\om;\la) \bigr)
- \cdots}{p'(x,|\xi|(x,\om;\la)\om)}.
\eqno (4.4)
$$

Write:
$$
p(x, |\xi|\om) \sim p_d(x, \om) |\xi|^d + p_{d-1}(x, \om) |\xi|^{d-1}
+ p_{d-2}(x, \om) |\xi|^{d-2} \cdots.
$$
To solve $p(x, |\xi|\om) = \la$ amounts to a series reversion.

In order to see how that is done, let us assume for a short while that
$H$ is a first-order operator with constant coefficients ---for
instance, the absolute value of the Dirac operator on $\R^n$. We then
expect
$$
|\xi|(x,\om;\la) \sim \frac{1}{p_1(\om)}\,\la
- \frac{p_0(\om)}{p_1(\om)} - p_{-1}(\om)\,\la^{-1} + \cdots.
$$
Integration over $|\om| = 1$ gives
$$
d_H(x,x;\la) \sim (2\pi)^{-n} \bigl(a_0\,\la^{n-1} + a_1 \,\la^{n-2}
+ a_2\,\la^{n-3} + \cdots \bigr) \ala C,
$$
where, clearly, $a_0 = \wres H^{-n}$.

To compute $a_1,a_2,\dots$ we can as well assume that the development
of $p$ is analytic as $|\xi| \to \infty$. Let
$\psi(z) := z^{n-1}/p'(z)$, so that
$$
a_0 \la^{n-1} + a_1(x) \la^{n-2} + a_2(x) \la^{n-3} + \cdots
\sim \int_{\Sf^{n-1}} \psi(|\xi|(x,\om;\la)) \,d\om.
$$
If $\Ga$ is a circle containing $|\xi|(x,\om;\la)$, wound once around
$\infty$, we have the Cauchy integral:
$$
\eqalign{
\psi(|\xi|(x,\om;\la))
&= \psi(p^{-1}(\la))
= - \frac{1}{2\pi i} \oint_\Ga \frac{\psi(z) p'(z)\,dz}{p(z) - \la}
\cr
&= \frac{1}{2\pi i} \oint_{\Ga^{-1}}
\frac{\psi(\zeta^{-1}) p'(\zeta^{-1})\,d\zeta}
{\zeta^2 (p(\zeta^{-1}) - \la)}
= \frac{1}{2\pi i} \oint_{\Ga^{-1}}
\frac{d\zeta}{\zeta^{n+1} (p(\zeta^{-1}) - \la)}.
\cr}
$$
Thus $a_j(x) = \int_{\Sf^{n-1}} c_j(x,\om) \,d\om$, where
$$
\eqalign{
c_j(\om) &= \frac{1}{2\pi i} \oint_{|s|=\eps}
s^{n-j-2} \psi(p^{-1}(1/s)) \,ds
\cr
&= \frac{1}{(2\pi i)^2} \oint_{|s|=\eps} s^{n-j-2} \,ds
\oint_{\Ga^{-1}} \frac{d\zeta}{\zeta^{n+1}\,(p(1/\zeta) - 1/s)} \cr
&= \frac{1}{(2\pi i)^2} \oint_{\Ga^{-1}} \frac{d\zeta}{\zeta^{n+1}}
\oint_{|s|=\eps} \frac{s^{n-j-1} \,ds}{s\,p(1/\zeta) - 1} \cr
&= \frac{1}{2\pi i} \oint_{\Ga^{-1}}
\frac{d\zeta}{\zeta^{n+1} p(1/\zeta)^{n-j}}, \cr}
$$
which is the coefficient of $\zeta^n$ in the expansion of
$p(1/\zeta)^{j-n}$. Integrating over $|\om| = 1$ yields thus
$$
a_j = \wres H^{j-n},
$$
so, finally:
$$
d_H(x,x;\la) \sim \frac{1}{(2\pi)^n}
(\wres H^{-n}\la^{n-1} + \wres H^{-n+1}\la^{n-2} + \cdots) \ala C,
$$
where the densities of Wodzicki residues are constant for a
constant-coefficient operator. It is amusing that we have arrived at a
version of the classical Lagrange--B\"urmann expansion
\cite{\Lagrange}, with Wodzicki residues in the place of ordinary
residues.

Notice that $a_n = 0$. This is a very simple ``vanishing theorem''
(see for instance~\cite{\BransonG}).

Returning to the general case, if $H$ is a positive pseudodifferential
operator of order~$d$, then $A := H^{1/d}$ is a positive
pseudodifferential operator of first order. Setting $\mu = \la^{1/d}$,
we have
$$
\dl(\la - H) = \dl(\mu^d - A^d) = \frac{\dl(\mu - A)}{d\mu^{d-1}}
= \frac{\dl(\la^{1/d} - H^{1/d})}{d\la^{(d-1)/d}}.
$$
and so
$$
d_H(x,x;\la) \sim \frac{1}{d\,(2\pi)^n}
\bigl( a_0(x) \la^{(n-d)/d} + a_1(x) \la^{(n-d-1)/d}
+ a_2(x) \la^{(n-d-2)/d} + \cdots \bigr) \ala C.
\eqno (4.5)
$$
Clearly, $a_0 = \wres H^{-n/d}$. Now, the order of $q_2$ is at most
$2d - 1$, therefore its higher order contribution to this development
is in principle to $a_1$; the order of $q_3$ is at most $3d - 2$, so
it contributes to $a_2$ at the earliest, and so on.

Formula~\eq{4.5}, obtained through fairly elementary manipulations, is
the main result of this section. To illustrate its power, we show how
to reap from it a rich harvest of classical results (with a little
extra effort).

\proclaim Corollary 4.1.
\textrm{(Connes' trace theorem)}\enspace
For positive elliptic pseudodifferential operators of order~$-n$ on a
compact $n$-dimensional manifold, the Dixmier trace and the Wodzicki
residue are proportional:
$$
\Dtr H = \frac{1}{n\,(2\pi)^n} \Wres H.
$$

\Proof.
Let $H$ be of order $d = -n$ in~\eq{4.5}. We get
$$
d_H(x,x;\la)
\sim -\frac{1}{n\,(2\pi)^n} \wres H(x) \,\la^{-2} + \cdots \ala C.
$$
Assume the manifold is compact. We then know that $H$ is a compact
operator. Now, heuristically the argument goes as follows:
$N'_H(\la) \sim -\la^{-2}$, ergo $N_H(\la) \sim \la^{-1}$, ergo
$\la_l(H) \sim l^{-1}$. A Tauberian argument can be used at this
point~\cite{\Sirius} to ensure that the second asymptotic estimate is
valid without the Ces\`aro condition; and then the result follows. But
this is by no means necessary. One can steal a look at Section 6
and, by approaching step functions by elements of $\SS$, prove in an
elementary way that for any given $\eps > 0$ there is $l(\eps)$ such
that
$$
\frac{C(1 - \eps)}{l(\eps)} < \la_l(H) < \frac{C(1 + \eps)}{l(\eps)},
$$
where $C = n^{-1}(2\pi)^{-n} \Wres H$.
\endproof

On a noncompact spin manifold, consider now the Dirac operator on the
space of spinors $L^2(S)$. The noncommutative integral of $|D|^{-n}$
does not exist. However, if $\int a(x) \,dx$ is defined, it is
computable by a noncommutative integral:

\proclaim Theorem 4.2.
Let $a$ be an integrable function with respect to the volume form
on~$M$. Then
$$
C_n \int_M a(x) \,dx = \frac{1}{n\,(2\pi)^n} \Wres(a|D|^{-n}),
$$
where on the right hand side $A$ is seen as a multiplication operator
on $L^2(S)$. The constants are $C_{2k} = (2\pi)^{-k}/k!$ and
$C_{2k+1} = \pi^{-k-1}/(2k + 1)!!$

\Proof.
That follows from Theorem 5.3 of~\cite{\Sirius} if $a$ is a smooth
function with compact support. For $a$ positive and integrable, use
monotone convergence on both sides; the general case follows at once.
\endproof

The former is a small step in the direction of a theory of $K$-cycles
(or ``spectral triples'', as they are nowadays called) over noncompact
manifolds.

\proclaim Corollary 4.3.
\textrm{(Weyl's estimate)}\enspace
Let $N_H(\la)$ denote the counting function of $H$, a Laplacian on a
compact manifold or bounded region $M$  acting on scalar functions.
Then
$$
N_H(\la) \sim \frac{\Om_n \vol M}{n(2\pi)^n}\la^{n/2},
$$
where $\Om_n$ is the surface area of the unit ball in $\R^n$.

\Proof.
The same type of arguments as in Corollary~4.1 work. Indeed, this
estimate is a corollary of it~\cite{\Sirius}.
\endproof

Next consider Schr\"odinger operators $-\Delta + V(x)$, with symbol
$p(x,\xi) = |\xi|^2 + V(x)$. We can take a slightly different tack and
solve the equation $p(x,\xi) = \la$ by
$|\xi| = \sqrt{(\la - V(x))_+}$.

\proclaim Corollary 4.4.
\textrm{(The correspondence principle)}\enspace
For Schr\"odinger operators:
$$
N_H(\la) \sim \frac{\Om_n}{n(2\pi)^n} \int(\la - V(x))_+^{n/2}\,dx.
\eqno\qed
$$

See~\cite{\Gurarie}, for instance, for the reasons for the terminology.

A word of caution is in order here. The development~\eq{4.5} cannot be
integrated term by term in general. Consider, for instance, the
harmonic oscillator hamiltonian $H = \thalf(- d^2/dx^2 + x^2)$
on~$\R$: according to the theory developed here, its spectral density
behaves as $1/\sqrt{\la}$. If $\psi_n$, $n \in \N$ denote the
normalized wavefunctions, then indeed, like in Fourier series theory,
$$
\sum_{n+\shalf\leq\la} \psi_n^2(x) \sim \frac{\sqrt{\la}}{\pi}
$$
is true and can be independently checked. But $\wres H^{-1/2}$ is not
integrable over the real line, so one cannot conclude that $N_H(\la)$
behaves as $\sqrt{\la}$. Actually, as we saw in Section~2,
$\sum_{n=0}^\infty \dl(\la - (n + \thalf))
= H(\la) + o(\la^{-\infty}) \ala C$, so
$N_H(\la) = \la H(\la) + o(\la^{-\infty}) \ala C$. Now, Corollary~4.2
applies, so we have
$$
N_H(\la) \sim \frac{2}{2\pi}
\int_{-\sqrt{2\la}}^{\sqrt{2\la}} \sqrt{2\la - x^2} \,dx = \la H(\la)
$$
precisely as it should. (See the discussion in~\cite{\Gianni}.)

Consider $n$-dimensional Schr\"odinger operators with (continuous)
homogeneous potentials $V(x) \geq 0$, $V(ax) = t^a V(x)$. The previous
formula gives
$$
N_H(\la) \propto \la^{n/2+n/a}\int_{\Sf^{n-1}}V(x)^{-n/a}\,dx.
$$
and this means that if the cone $\set{x \in \R^n : V(x) = 0}$ is too
big, in the counting number estimate we are heading for
trouble~\cite{\Solomyak}. But the ``nonstandard asymptotics'' that
might then intervene do not detract from the validity of the
nonintegrated formula~\eq{4.5}.

\smallskip

In the remainder of the section, we focus on the computation of
spectral densities for Laplacians. Nothing essential is won or lost by
considering general vector bundles, so we work on scalars. The more
general Laplacian operator on a Riemannian manifold is (minus) the
Laplace--Beltrami operator $\Delta$ plus potential vector and scalar
potential terms, with symbol
$$
\eqalign{
p(x,\xi)
&= -g^{ij}(x)\bigl(\xi_i\xi_j + (i\Ga_{ij}^k(x) \xi_k + 2A_i(x) \xi_j)
\cr
&\qquad + (A_i(x)A_j(x) + i(\Ga_{ij}^k(x)A_k(x)
- \del_i A_j(x))) \bigr) + V(x)
\cr
& =: -g^{ij}(x)\xi_i\xi_j + B^i(x)\xi_i + C(x).
\cr}
$$
Formula~\eq{4.5} would seem to give for this case:
$$
d_H(x,x;\la) \sim \frac{1}{2\,(2\pi)^n}
\bigl( a_0(x) \la^{(n-2)/2} + a_1(x) \la^{(n-3)/2}
+ a_2(x) \la^{(n-4)/2} + \cdots \bigr) \ala C.
$$
In fact, it will be seen in a moment that $a_1 = a_3 = \cdots = 0$.
Also we know already that $a_0(x) = \Wres \Delta^{-n/2} = \Om_n$. Our
task is to compute the next coefficients; it is a rather exhausting
one, whose results can be inferred from the extensive work already
carried out~\cite{\Gilkey} on heat kernel expansions (see Section~6),
so we will limit ourselves to the computation of $a_2(x)$ to
illustrate the relative simplicity of our approach.

Let $n \geq 3$. Write $a$ for $g^{ij}(x)\om_i\om_j$, then $b$ for
$B^i(x)\om_i$ and $c$ for $C(x)$. Our method calls for solving for the
positive root of $a|\xi|^2 + b |\xi| + (c - \la) = 0$ and substituting
this in $|\xi|^{n-1}/(2a|\xi| + b)$. In diminishing powers of $\la$,
we obtain for the latter the development:
$$
\frac{1}{2a^{n/2}} \biggl(\la^{(n-2)/2}
- \frac{(n-1)b}{2a^{1/2}} \la^{(n-3)/2}
+ \Bigl(\frac{n(n-2)b^2}{8a} - \frac{(n-2)c}{2}\Bigr) \la^{(n-4)/2}
+ \cdots \biggr).
\eqno (4.6)
$$
One sees that odd-numbered terms in this expansion contain odd powers
of $\om$ and thus give vanishing contributions, after the integration
on the cosphere. Also, the contribution of the $q_2$ term in~\eq{4.2}
will start at order $\thalf n - 2$ in $\la$, the contribution of $q_3$
will start at order $\thalf n - 3$ and so on: the terms in the
asymptotic expansion of the density kernels of Laplacian operators
differ by powers of $\la$, not of $\sqrt{\la}$, as one would expect on
general grounds.

It is convenient now to use geodesic coordinates at each point; this
is justified by the nature of the result. In these coordinates
$\Ga^k_{ij}(x_0) = 0$ and we have the Taylor expansion
$$
g_{ij}(x) \sim \dl_{ij}
+ \frac{1}{3} R_{iklj}(x_0) \,(x-x_0)^k (x-x_0)^l
+ \sum_{|\a|\geq3} \del^\a g(x_0) \frac{(x-x_0)^\a}{\a!}
\as x \to x_0,
$$
where $R_{iklj}$ denotes the Riemann curvature tensor. Recall that the
Ricci tensor is given by $R_{kj} := \sum_l R^l_{klj}$ and the scalar
curvature by $R := \sum_{kj} g^{kj}R_{kj}$.

{}From~\eq{4.3} one obtains for $q_2(x_0, \xi)$
$$
\frac{1}{2} \sum_{|\a|>0}\frac{i^{-|\a|}}{\a!}
\del_\xi^\a	\bigl(-g^{ij}(x_0)\xi_i\xi_j + B^i(x_0)\xi_i\bigr)
\, \del_x^\a \bigr|_{x=x_0}
\bigl(-g^{ij}(x)\xi_i\xi_j + B^i(x)\xi_i + C(x)\bigr).
\eqno (4.7)
$$
Let us take for a moment $A_i = 0$. Then in geodesic coordinates
$B^i(x_0) = 0$ and it is not hard to see that the only surviving term
in~\eq{4.7} is equal to $\frac{1}{3} R_{kj}(x_0) \xi^k\xi^j$. Also
$b = 0$ in~\eq{4.6}. So, in view of~\eq{4.4} we are left with two
terms at order $\la^{(n-4)/2}$, to wit:
$$
-\int_{\Sf^{n-1}} \! d\om\, \frac{(n-2)C(x_0)}{2} \la^{(n-4)/2}
$$
that comes from the third term in~\eq{4.6}, and the first order
contribution of
$$
\int_{\Sf^{n-1}} \! d\om\, \frac{\frac{\del^2}{\del\la^2}
\bigl( q_2(x,|\xi|(x,\om;\la)\om)|\xi|^{n-1} (x,\om;\la) \bigr)}
{p'(x,|\xi|(x,\om;\la)\om)}.
$$
In effect, $q_2$ contributes here a factor of order $\la$, so the
second derivative in the previous formula gives rise to a term of
order $\la^{(n-4)/2}$ also. To finish the computation, we use
$$
\int_{\Sf^{n-1}} \! d\om\, A_{ij}\om^i\om^j
= \frac{\Om_n}{n} g^{ij} A_{ij},
$$
to get
$$
a_2(x_0)
= \frac{(n - 2)\Om_n}{2} \bigl(\tfrac{1}{6} R(x_0) - C(x_0)\bigr).
\eqno (4.8)
$$
Notice that for a pure Laplace--Beltrami operator, the contribution to
$a_2$, when computed in geodesic coordinates, comes exclusively
through the $q_2$ term.

It remains to convince ourselves that vector potentials give no
contribution at this stage. On one hand, the $c$ term in~\eq{4.6}
would contribute now the extra terms
$$
-\frac{(n - 2)\Om_n}{2} (A^jA_j + i\del_jA^j).
$$
On the other, the term in $b^2$ in the same formula would contribute a
term of the form $\thalf(n - 2)\Om_n A^jA_j$, and in the computation
of $q_2$ there appears now a term $(2i/n)\,\del_jB^j$ that contributes
$\tihalf (n - 2)\Om_n\,\del_jA^j$ and thereby cancels the rest.
Therefore~\eq{4.8} stands also in that case.

Actually the coefficients of the Ces\`aro asymptotic expansion of
$d(x,x;\la)$ are all (local densities of) Wodzicki residues for $n$
odd: $a_{2k}(x) = \wres\Delta^{-n/2+k}(x)$, for $k \in \N$. For $n$
even we have $a_{2k} = \wres\Delta^{-n/2+k}$ only as long as
$-n/2 + k < 0$ (the Wodzicki residues of nonnegative powers of a
differential operator being of course zero); the following
coefficients for the parametric expansion are, in our terminology of
Section~2 (further explained in the next two sections), not
``residues'' but ``moments''. Note that for $n = 2$, the coefficient
$a_2$ is already a ``moment'' and cannot be computed by a Ces\`aro
development. This strikingly different behaviour of the
odd-dimensional and the even-dimensional cases is concealed in the
uniformity of the usual heat kernel method, but it reflects itself in
the fact that the corresponding zeta functions have an infinite number
of poles, corresponding to the residues, in the odd-dimensional case;
and a finite number in the even-dimensional case. One
has~\cite{\Wodzicki}:
$$
\Res_{s=n/2-k} \zeta_H(s) = \thalf \Wres H^{k-n/2},
$$
where
$$
\zeta_H(s) = \int_M \bigl< d_H(x,x;\la), \la^{-s} \bigr>_\la \,dx
\qquad (\Re\,s \gg 0)
$$
is the kernel of the zeta operator~\eq{3.1}. A direct, ``elementary''
proof of the essential identity between Wodzicki residues and residues
of the poles of the zeta functions is obviously in the cards, but we
will not go further afield here. For a nontrivial use of the
noncommutative residue in zeta function theory, have a look
at~\cite{\ElizaldeVZ}.

\section 5. Ces\`aro developments of counting functions

We consider here operators on compact manifolds without boundary and
look at the behaviour of the counting function
$$
N(\la) := \sum_{\la_l\leq\la} 1.
$$

In order to refresh our intuition, we shall follow a deliberately
na\"{\i}ve approach and temporarily forget some of what we learned at
the end of last section. Envisage first the scalar Laplacian on $\T^2$
with the flat metric; then the counting function is given by the
following table:
$$
\vbox{\offinterlineskip
\halign{\vrule#&\strut\hfil$#$\hskip 3pt\vrule\hskip 1pt
&&\vrule\hskip 3pt$#$\hskip 3pt\cr
\omit&\la & 0 & 1 & 2 & \hfil 4 & \hfil 5 & \hfil 8 & \hfil 9 & 10
& 13 & 16 & 17 & 18 & 20 & 25 & 26 & \cdots\cr height
2pt&\omit\hfil\hskip 3pt\vrule\hskip 1pt
&\omit\vrule\hfil&\omit\vrule\hfil&\omit\vrule\hfil
&\omit\vrule\hfil&\omit\vrule\hfil&\omit\vrule\hfil
&\omit\vrule\hfil&\omit\vrule\hfil&\omit\vrule\hfil
&\omit\vrule\hfil&\omit\vrule\hfil&\omit\vrule\hfil
&\omit\vrule\hfil&\omit\vrule\hfil&\omit\vrule\hfil
&\omit\vrule\hfil\cr
\noalign{\hrule}
height 2pt&\omit\hfil\hskip 3pt\vrule\hskip 1pt
&\omit\vrule\hfil&\omit\vrule\hfil&\omit\vrule\hfil
&\omit\vrule\hfil&\omit\vrule\hfil&\omit\vrule\hfil
&\omit\vrule\hfil&\omit\vrule\hfil&\omit\vrule\hfil
&\omit\vrule\hfil&\omit\vrule\hfil&\omit\vrule\hfil
&\omit\vrule\hfil&\omit\vrule\hfil&\omit\vrule\hfil
&\omit\vrule\hfil\cr
\omit&N(\la^+) & 1 & 5 & 9 & 13 & 21 & 25 & 29 & 37
& 45 & 49 & 57 & 61 & 69 & 81 & 89 & \cdots\cr}}
$$
No doubt, $N(\la) \sim \pi\la$ is a reasonable first approximation;
but it is also plain that the remainder undergoes wild oscillations.
The precise determination of this remainder is a difficult problem,
not unlike the problem of determining the next-to-main term in the
asymptotic development of prime numbers.

An even simpler and more telling example is provided by the
eigenvalues $\la_l$ of the Laplacian on the $n$-dimensional sphere.
They are given by
$$
\la_l = l(l + n - 1)
\hskip-1em \sepword{with respective multiplicities}
m_l = {l + n \choose n} - {l + n - 2 \choose n},
\eqno (5.1)
$$
for $l \in \N$. For example, if $n = 2$, the eigenvalues are $l(l+1)$
and the multiplicities are $(2l + 1)$. The leading term is
$$
N(\la) \sim \frac{2}{n!}\, \la^{n/2} \as \la \to \infty.
$$
On the other hand, asymptotically:
$$
N(\la^+) - N(\la^-) \sim \frac{2\,l^{n-1}}{(n-1)!},
$$
and so
$$
\la^{(1-n)/2} \bigl(N(\la^+) - N(\la^-)\bigr) \sim \frac{2}{(n-1)!}.
$$
Plainly, we cannot find an asymptotic formula for $N(\la)$ with error
term $o(\la^{(n-1)/2})$ and continuous main term. The example is taken
from H\"ormander's work~\cite{\HormanderS, \HormanderIII}.

The foregoing is a ``Gibbs phenomenon'' related to the lack of
smoothness of the characteristic function. The problem is ``solved''
if one is prepared to look at the expansions in the Ces\`aro sense.
The fact that higher order terms in the asymptotic expansion of the
eigenvalues of the Laplacian were to be understood in an averaged
sense was pointed out by Brownell~\cite{\Brownell} many years ago.

Going back to tori, consider the distribution of nonvanishing
eigenvalues $\{\la_l\}_{l=1}^\infty$ of the scalar Laplacian on an
$n$-dimensional torus $\T^n$, with the flat metric. The eigenfunctions
$\{\phi_l\}_{l=1}^\infty$ can be seen as nonzero smooth functions in
$\R^n$ that satisfy
$$
\Delta \phi_l + \la_l \,\phi_l = 0
$$
and the periodicity conditions
$$
\phi_l(x_1 + 2k_1\pi, \dots, x_n + 2k_n\pi) = \phi_l(x_1,\dots,x_n),
$$
where the girths of the torus are taken to be $2\pi$ in all
directions.

Those eigenvalues are given by $\la_k = k_1^2 +\cdots+ k_n^2$ for
$k = (k_1,\dots,k_n) \in \Z^n$, with corresponding eigenfunctions
$\phi_k(x_1,\dots,x_n) = e^{ik\cdot x}$. Thus the $\la_l$ are the
nonnegative integers $q_l$ that can be written as a sum of $n$
squares. The multiplicity of each such value is the number of integral
solutions of the Diophantine equation $q_l = k_1^2 +\cdots+ k_n^2$. We
wish to compute the terms in the parametric and Ces\`aro developments
of $N(\la)$ next to leading Weyl term (which in fact for this problem
goes back to Gauss):
$$
N(\la) \sim \frac{\Om_n}{n} \,\la^{n/2} \as \la \to \infty.
$$

To do so, we start with the derivative $N'(\la)$; this is
nothing but $(2\pi)^n d(x, x; \la)$, but, as advertised, it is more
instructive to forget for a while the discussion in Section~4. We
have:
$$
N'(\la) = \sum_{l=1}^\infty \dl(\la - \la_l)
= \sum_{k\in\Z^n} \dl(\la - k_1^2 - \cdots - k_n^2).
$$
Let $\phi \in \D(\R)$, let $\sigma$ be a large real parameter and set
$\eps = 1/\sigma$, so that $\eps \downarrow 0$. Then
$$
\eqalign{
\bigl< N'(\sigma\la), \phi(\la) \bigr>_\la
&= \eps \bigl< N'(x), \phi(\eps\la) \bigr>_\la
= \eps \sum_{k\in\Z^n} \phi(\eps|k|^2)
\cr
&= \eps^{1-n/2} \int_{\R^n} \phi(|x|^2) \,dx + o(\eps^\infty)
\cr
&= \thalf \Om_n \,\eps^{1-n/2} \int_0^\infty r^{(n-2)/2} \phi(r) \,dr
+ o(\eps^\infty).
\cr}
$$
The third equality is just Lemma~2.10.

Hence, weakly:
$$
N'(\sigma\la)
= \thalf \Om_n \sigma^{-1+n/2} \la_+^{-1+n/2} + o(\sigma^{-\infty})
\as \sigma \to \infty,
$$
and upon integration
$$
N(\sigma\la)
= \frac{\Om_n}{n} \,\la_+^{n/2} \,\sigma^{n/2} + o(\sigma^{-\infty})
\as \sigma \to \infty.
$$
Observe that the constant of integration $\mu_0$ vanishes, as do all
the other moments.

Then Theorem 2.4 yields:
$$
N(\la) = \frac{\Om_n}{n} \,\la^{n/2} + o(\la^{-\infty}) \ala C
\as \la \to \infty.
$$
Hence the error term, although definitely not small in the ordinary
sense, is of rapid decay in the $(C)$ sense.

We turn to examine some cases of spheres. The derivative of the
counting function for $\Sf^2$ is
$N'(\la) = \sum_{l=0}^\infty (2l + 1) \,\dl(\la - l(l+1))$. To deal
with this case, we need a heavier gun than Lemmata 2.9--2.11. This is
provided by:

\proclaim Lemma 5.1.
Let $f \in \K'(\R^n)$, so that it satisfies the moment asymptotic
expansion. If $p$ is an elliptic
polynomial and $\phi \in \SS$, then
$$
\bigl< f(x), \phi(t p(x)) \bigr>
\sim \sum_{m=0}^\infty \frac{\<f(x),p(x)^m>\,\phi^{(m)}(0)}{m!} \,t^m
\as t \to 0.
$$

\Proof.
The proof consists in showing that the Taylor expansion
$$
\phi(t p(x))
= \sum_{m=0}^N \frac{\phi^{(m)}(0) p(x)^m}{m!} \,t^m + O(t^{N+1})
$$
holds not only pointwise, but also in the topology of~$\K(\R^n)$.
\endproof

Consider now the distribution
$$
f(\la)
:= (2\la + 1) \biggl( \sum_{l=1}^\infty \dl(\la - l) - H(\la) \biggr),
$$
that lies in $\K'$. Notice that
$$
\eqalign{
\bigl< f(\la), \phi(t(\la^2 + \la)) \bigr>
&= \sum_{l=1}^\infty (2l + 1)\,\phi(t(l^2 + l))
- \int_0^\infty (2\la + 1) \phi(t(\la^2 + \la)) \,d\la
\cr
&= \sum_{l=1}^\infty (2l + 1)\,\phi(t(l^2 + l))
- \int_0^\infty \phi(t\mu) \,d\mu.
\cr}
$$
{}From Lemma~5.1 we conclude that, for $\phi \in \SS$,
$$
\eqalign{
\bigl< N'(\la), \phi(t\la) \bigr>
&= \sum_{l=0}^\infty (2l + 1)\,\phi(t(l^2 + l))
\cr
&\sim \int_0^\infty \phi(t\mu) \,d\mu + \phi(0) + \sum_{j=0}^\infty
\frac{\<f(\la),(\la^2 + \la)^j>\,\phi^{(j)}(0)}{j!} \,t^j
\as t \downarrow 0.
\cr}
$$

The parametric expansion of $N'(\la)$ is thus
$$
N'(\la/t) \sim H(\la) + \dl(\la)t
+ \sum_{j=0}^\infty \frac{(-1)^j\mu_j\,\dl^{(j)}(\la)}{j!} \,t^{j+1}
\as t \downarrow 0,
$$
where the ``generalized moments'' $\mu_j$ are given by
$$
\mu_j = \bigl< f(\la), (\la^2 + \la)^j \bigr>
= \sum_{l=1}^\infty (2l + 1) (l^2 + l)^j
- \int_0^\infty (2\la + 1) (\la^2 + \la)^j \,d\la  \ala C.
$$
It follows that $N'(\la) \sim H(\la) + o(\la^{-\infty}) \ala C$ as
$\la \to \infty$.

In view of our gymnastics with Riemann's zeta function in Section~2,
the computation of the $\mu_j$ presents no difficulties. We obtain
$$
\eqalign{\mu_0 &= 2\zeta(-1) + \zeta(0) = -\frac{2}{3}, \cr
\mu_1 &= 2\zeta(-3) + \zeta(-1) = -\frac{1}{15}, \cr} \quad
\eqalign{\mu_2 &= 2\zeta(-5) + 4\zeta(-3) = \frac{8}{315}, \cr
\mu_3 &= 2\zeta(-7) + 9\zeta(-5) + \zeta(-3) = -\frac{2}{105},}
$$
and so on. On integrating, we get
$$
N(\la/t) \sim \frac{\la}{t} H(\la) + \frac{1}{3} H(\la)
+ \frac{1}{15} \dl(\la) \,t + \frac{4}{315} \dl'(\la) \,t^2 + \cdots
\as t \downarrow 0,
\eqno (5.2)
$$
and $N(\la) \sim \la\,H(\la) + \frac{1}{3}\,H(\la)
+ o(\la^{-\infty}) \ala C$. Note that the $\la^0$th order term in the
Ces\'aro development for $N(\la)$ comes from the first moment. The
curvature of a sphere $\Sf^n$ is given by $R = n(n-1)$, so the second
term in the development is precisely what we had expected.

\smallskip

We look now at the derivative of the counting function for the
Laplace--Beltrami operator on $\Sf^3$. It is slightly simpler to
consider the operator $1 - \Delta$, for which we have, according
to~\eq{5.1}:
$N'(\la) = \sum_{l=0}^\infty (l + 1)^2 \,\dl(\la - (l+1)^2)$.

Consider the distribution
$$
f(\la)
:= (\la+1)^2 \biggl(\sum_{l=0}^\infty \dl(\la - l) - H(\la+1)\biggr),
$$
lying in $\K'$. We have:
$$
\bigl< f(\la), \phi(t(\la + 1)^2) \bigr>
= \sum_{l=0}^\infty (l + 1)^2 \phi(t(l + 1)^2)
- \int_{-1}^\infty (\la + 1)^2 \phi(t(\la + 1)^2) \,d\la.
$$
One sees that the moments all cancel:
$\<f(\la), (\la + 1)^{2j}> = \zeta(-2j - 2) = 0$, for $j \in \N$.
Therefore we get simply
$$
\bigl< N'(\la), \phi(t\la) \bigr>
\sim \frac{1}{2t^{3/2}} \int_0^\infty \phi(u) \sqrt{u} \,du
\as t \downarrow 0,
$$
and thus in this case we collect just the Weyl term
$$
N(\la) \sim \frac{\la^{3/2}H(\la)}{3} \ala C  \as \la \to \infty.
\eqno (5.3)
$$

We may reflect now that the counting number for these Laplacians on
$\Sf^2$, $\Sf^3$ behave in the expected way for even and odd
dimensional cases, respectively. For a generalized Laplacian which is
the square of a Dirac operator the qualitative picture is the same. In
particular, the Chamseddine--Connes expansion corresponds to $n = 4$,
whereupon the counting functional behaves in much the same way as the
one for~$\Sf^2$. Therefore, formal application of the
Chamseddine--Connes Ansatz to the characteristic function of the
spectrum, as done in~\cite{\CarminatiIKS, \IochumKS} misses the terms
involving $\dl$ and its derivatives ---whose physical meaning, if any,
is unclear to us.

\section 6. Spectral density and the heat kernel

Now we tackle the issue of the small-$t$ behaviour of the Green
functions associated to an elliptic pseudodifferential operator~$H$.
These are the integral kernels of operator-valued functions of $H$, of
the form
$$
G(t,x,y) = \bigl< d_H(x,y;\la), g(t\la) \bigr>_\la
$$
where $g$, as already advertised, will in this section belong (or can
be extended) to the Schwartz space $\SS$ (i.e., we deal with the
standard theory as opposed to the framework sketched at the end of
Section 3). The basic question is whether $G(t,x,y)$ has an asymptotic
expansion as $t \downarrow 0$. In effect, we shall see immediately how
to obtain from the $(C)$ asymptotic expansion for the spectral density
an \textit{ordinary} asymptotic expansion for Green functions.

The emphasis in recent years has been on Abelian type expansions, the
so-called heat kernel techniques~\cite{\Gilkey}. It is common folklore
that Ces\`aro summability implies Abel summability, but not
conversely. As we just claimed, one can go from the Ces\`aro expansion
to the heat kernel expansion. The reverse implication does not work
quite the same. If we know the coefficients of the heat kernel
expansion and we independently know that a Ces\`aro type expansion for
the spectral density exists, we can infer the coefficients of the
latter from the former. But it may happen that the formal
Abel--Laplace type expansion does not say anything about the ``true''
asymptotic development.

For instance, if $f(\la) := \sin \la \,e^{\sqrt\la}$ for $\la > 0$,
then $\lim_{\la\to\infty} f(\la) \ala C$ does not exist, since no
primitive of~$f$ can have polynomial order in~$\la$. Even so, one can
show that $k(t) = \<f(\la), e^{-t\la}>$ has a Laplace expansion
$k(t) \sim a_{-1}t^{-1} + a_0 + a_1 t + \cdots$ as $t \downarrow 0$,
that is, $\lim_{\la\to\infty} f(\la) = a_{-1} \ala A$. To get an
example of a bounded function with this behaviour, one uses the fact
that $f_m(\la) = \sin \la^{1/m}$ obeys
$\lim_{\la\to\infty} f_m(\la) = 0 \ala{C,N}$ only for $N > m$,
together with Baire's theorem, to construct a bounded function
$f(\la) = \sum_{k\geq 1} 2^{-k} f_{m_k}(\la)$ that does not have a
Ces\`aro limit as $\la \to \infty$, but for which $f(\la) \to 0$ in
the Abel sense.

In order to relate our Ces\`aro asymptotic expansions with heat kernel
developments, we need to examine expansions of distributions $f(\la)$
that may contain nonintegral powers of~$\la$. Suppose that
$\{\a_k\}_{k\geq 1}$ is a decreasing sequence of real numbers, not
including negative integers, and suppose further that $f \in \SS'$,
supported in $[0,\infty)$, has the Ces\`aro asymptotic expansion
$$
f(\la) \sim \sum_{k\geq 1} c_k \la^{\a_k}
+ \sum_{j\geq 1} b_j \la^{-j} \ala C \as \la \to \infty.
$$
It follows from Theorem~32 of~\cite{\EstradaKbook} and from
Theorem~2.5 that $f$ has the following parametric development:
$$
f(\sigma\la) \sim \sum_{k\geq 1} c_k (\sigma\la_+)^{\a_k}
+ \sum_{j\geq 1} b_j \Pf((\sigma\la)^{-j} H(\la))
+ \sum_{m\geq 0} \frac{(-1)^m\mu_m\,\dl^{(m)}(\la)}{m!\,\sigma^{m+1}}
\eqno (6.1)
$$
as $\sigma \to \infty$, where the ``generalized moments'' $\mu_m$ are
given by
$$
\mu_m = \bigl< f(x) - \sum_{k\geq1} c_k x_+^{\a_k}
- \sum_{j\geq1} b_j \Pf(x^{-j}H(x)),\,x^m \bigr>
\eqno (6.2)
$$
and where $\Pf$ denotes a ``pseudofunction'' \cite{\EstradaKfp}
obtained by taking the Hadamard finite part, that is:
$\<\Pf(h(x)),g(x)> := \Fp\int_0^\infty h(x) g(x)\,dx$ if
$\supp h \subseteq [0,\infty)$. In particular,
$$
\eqalign{
&\bigl< \Pf(x^{-j}H(x)), g(x) \bigr>
= \Fp\int_0^\infty \frac{g(x)}{x^j} \,dx \cr
&\qquad = \int_1^\infty \frac{g(x)}{x^j} \,dx
+ \int_0^1 \frac{1}{x^j} \biggl( g(x)
- \sum_{k=0}^{j-1} \frac{g^{(k)}(0)}{k!} x^k \biggr) \,dx -
\sum_{k=0}^{j-2} \frac{g^{(k)}(0)}{k!(j-k-1)}. \cr}
\eqno (6.3)
$$
Notice that taking the finite part involves dropping a logarithmic
term proportional to $g^{(j-1)}(0)$. This has the consequence that
$\Pf(x^{-j}H(x))$ fails to be homogeneous of degree~$-j$ by a
logarithmic term; indeed,
$$
\Pf((\sigma\la)^{-j} H(\sigma\la))
= \sigma^{-j} \Pf(\la^{-j} H(\la))
+ \frac{(-1)^j \dl^{(j-1)}(\la) \,\log\sigma}{(j-1)!\,\sigma^j}.
$$
Consequently,
$$
\eqalignno{
\bigl< f(\la), g(t\la) \bigr>_\la
&\sim \sum_{k\geq 1} c_k \,t^{-\a_k-1}
\Fp\int_0^\infty \la^{\a_k} g(\la) \,d\la \cr
&\qquad + \sum_{j\geq 1} b_j\,t^j \biggl(
\Fp\int_0^\infty \frac{g(\la)}{\la^j} \,d\la -
\frac{g^{(j-1)}(0)}{(j-1)!} \log t \biggr) \cr
&\qquad + \sum_{m\geq 0} \frac{\mu_m \,g^{(m)}(0)}{m!}\,t^m.
& (6.4) \cr}
$$

The heat kernel development may be recovered by taking
$g(\la) = e^{-\la}$ for $\la \geq 0$. In that case, $d\a_k$ is
integral, $\Fp\int_0^\infty \la^{\a_k}g(\la) \,d\la = \Ga(\a_k+1)$ and
$g^{(j-1)}(0) = (-1)^{(j-1)}$. {}From this it is clear that the heat
kernel of a pseudodifferential operator may generally contain
logarithmic terms. Indeed, by harking back to~\eq{4.5}, on using
\eq{6.4} we prove:

\proclaim Corollary 6.1.
The general form of the (coincidence limit of) the heat kernel for an
elliptic pseudodifferential operator of order $d$ on a compact
manifold $M$ of dimension $n$ is given by
$$
K(t,x,x) \sim \sum_{j-n\notin d\N_+} \!\! \ga_{j-n}(x) t^{(j-n)/d}
+ \sum_{j-n\in d\N_+} \!\! \b_{j-n}(x) t^{(j-n)/d} \log t
+ \sum_{r=1}^{\infty} r_m(x) t^m
$$
qs $t \downarrow 0$, where
$$
\ga_{j-n}(x) = \frac{\Ga((n - j)/d)}{d(2\pi)^n}\, a_j(x)
,
$$
and similarly for the other coefficients.

(See~\cite{\Grubb, Cor.~4.2.7}.)

Now suppose we know {\it a priori\/} that $f(\la)$ has a Ces\`aro
asymptotic expansion in falling powers of~$\la$, and that we also know
that $\Phi(t) := \bigl< f(\la), e^{-t\la} \bigr>_\la$ has an
asymptotic expansion as $t \downarrow 0$ without $\log t$ terms. Then
it follows that all $b_j = 0$ in~\eq{6.1}, i.e., there are no negative
integral exponents in the Ces\`aro development of~$f$, and
consequently the constants $\mu_m$ are the moments of~$f$.
Thus~\eq{6.4} simplifies to
$$
\Phi(t) \sim \sum_{k\geq 1} c_k \,\Ga(\a_k+1) \,t^{-\a_k-1}
+ \sum_{m\geq 0} \frac{(-1)^m \mu_m}{m!} \,t^m.
$$
This is precisely the case for a (generalized) Laplacian: if $n$ is
{\it odd}, only half-integer powers of $\la$ appear in the spectral
density and logarithmic terms in the heat kernel are thereby ruled
out. Notice that the Ces\`aro development for an odd dimensional
Laplacian need not terminate. For {\it even\/} dimensions, the term
$k = n/2$ is proportional to $\wres H^0 \la^{-1}$ and later terms are
proportional to $\wres H^r \la^{-r-1}$. However, since $H^r$ is a
differential operator, its local Wodzicki residue vanishes for
$r \in \N$, and the Ces\`aro development terminates at the $\la^0$
term. However, as we have seen, at this point the moments~\eq{6.2}
enter the picture.

\smallskip

It has become a habit to write the diagonal of the heat kernel for a
Laplacian in the form
$$
K(t,x,x) \sim (4\pi t)^{-n/2} \sum_{k=0}^\infty b_k(x,x) \,t^{k/2},
$$
where $n$ is the dimension of the manifold and $b_0(x,x) = 1$. We see
now that $b_k(x,x) = 0$ for $k$ odd, whereas
$$
b_{2k}(x,x)
= \frac{2^k a_{2k}(x)}{\Omega_n (n - 2)(n - 4)\dots (n - 2k)}
\qquad\hbox{for $k > 0$}.
$$
A similar formula holds off-diagonal. As we have noted, these
expansions are local in the sense that they do not distinguish between
a finite and an infinite region of $\R^n$, say. However, the smallness
of the terms after the first is not uniform near the boundary, and
hence the ``partition function''
$$
K(t) := \int_M K(t,x,x) \,dx
\sim (4\pi t)^{-n/2} \sum_{k=0}^\infty b_k \,t^{k/2},
\eqno (6.5)
$$
with $b_0 = \vol(M)$ for scalars, has an expansion with nontrivial
boundary terms in general, starting to contribute in
$b_2$~\cite{\McKeanS}.

As for the examples, the expansion \eq{6.5} for $\Sf^2$ was first
obtained as the partition function of a diatomic
molecule~\cite{\Mulholland} and is well known to physicists. On using
$\vol(\Sf^2) = 4\pi$, we read Mulholland's expansion directly by
looking at~\eq{5.2}:
$$
K_{\Sf^2}(t) \sim {1\over t} + \frac{1}{3} + \frac{1}{15}\,t
+ \frac{4}{315}\,t^2 + \cdots \as t \downarrow 0.
$$
As for the $SU(2)$ group manifold, from~\eq{5.3}, on using
$\vol(\Sf^3) = 2\pi^2$ and $e^{t\Delta} = e^t e^{-t(1-\Delta)}$, the
partition function is seen immediately to be
$$
K_{\Sf^3}(t) \sim {{\sqrt\pi}\over4t^{3/2}}e^t.
$$

\smallskip

We turn at last to the Chamseddine--Connes expansion. The theory of
Ces\`aro and parametric expansions justifies~\eq{1.1}, in the
following way. We work in dimension $n = 4$ and take $H = D^2$, a
generalized Laplacian, acting on a space of sections of a vector
bundle $E$, over a manifold without boundary. The kernel of its
spectral density satisfies
$$
d_{D^2}(x,x;\la) \sim \frac{\rk E}{16\pi^2} \la
+ {1 \over 32\pi^4} \wres D^{-2}(x) \ala C \as \la \to \infty.
$$
Integrating over~$M$ and using the formulas of this section with
$t = \La^{-2}$, we then get
$$
\eqalign{
\Tr \phi(D^2/\La^2) &\sim \frac{1}{(4\pi)^2}
\biggl( \rk E\, \La^4 \int_0^\infty \la \phi(\la) \,d\la
+ b_2 \La^2 \int_0^\infty \phi(\la) \,d\la
\cr
&\hskip5em + \sum_{m\geq 0} (-1)^m \phi^{(m)}(0)\,
b_{2m+4}(D^2) \,\La^{-2m} \biggr)  \as t \downarrow 0.
\cr}
$$
where $(-1)^m b_{2m+4}(D^2) = 16\pi^2 \mu_m(D^2)/m!$ are suitably
normalized, integrated moment terms of the spectral density of $D^2$.
Thus, we arrive at~\eq{1.1}.

We finally take stock of the status of the Chamseddine--Connes
development. If $\phi \in \SS$, then the development becomes a bona
fide asymptotic expansion. However, if one wishes to use (for
instance) the counting function $N_{D^2}(\la \leq \La^2)$, which does
not lie in~$\SS$, then the present formulae are not directly
applicable and one one must proceed like in Section 5; moreover the
expansion beyond the first piece is only valid in the Ces\`aro sense.
We close by noting that third piece of the Chamseddine--Connes
Lagrangian has interesting conformal properties; this is better
studied through the corresponding zeta function at the
origin~\cite{\Rosenberg}. That term is definitely not a Wodzicki
residue but a moment; whether this fact has any physical
significance is not easy to say.

\bigskip

\noindent{\bf Acknowledgments}

\medskip

Heartfelt thanks to S. A. Fulling for sharing his ideas with us prior
to the publication of~\cite{\EstradaF}. We wish to thank M.~Asorey,
E.~Elizalde, H.~Figueroa, D.~Kastler, F.~Lizzi, C.~P. Mart\'{\i}n,
A.~Rivero, T.~Sch\"ucker and J.~Sesma for fruitful discussions and
G.~Landi for a question that motivated the paragraph on harmonic
oscillators in Section~4. JMGB and JCV acknowledge support from the
Universidad de Costa Rica; JMGB also thanks the Departamento de
F\'{\i}sica Te\'orica de la Universidad de Zaragoza and JCV the Centre
de Physique Th\'eorique (CNRS--Luminy) for their hospitality.

\bigskip\bigskip

\leftline{\textbf{References}}

\bigskip
\frenchspacing

\refno\Ackermann.
Ackermann, T.:
A note on the Wodzicki residue.
J. Geom. Phys. {\bf 20}, 404--406 (1996)

\refno\Persephone.
Alvarez, E., Gracia-Bond\'{\i}a, J. M., Mart\'{\i}n, C. P.:
A renormalization group analysis of the NCG constraints
$m_{\rm top} = 2\,m_W$, $m_{\rm Higgs} = 3.14\,m_W$.
Phys. Lett. B {\bf 329}, 259--262 (1994)

\refno\BransonG.
Branson, T. P., Gilkey, P. B.:
Residues of the eta function for an operator of the Dirac type.
J. Func. Anal. {\bf 108}, 47--87 (1992)

\refno\Brownell.
Brownell, F. H.:
Extended asymptotic eigenvalue distributions for bounded domains in
$n$-space.
J. Math. Mech. {\bf 6}, 119--166 (1957)

\refno\Carleman.
Carleman, T.:
\"Uber die asymptotische Verteilung der Eigenwerte partieller
Diffe\-rentialgleichungen.
Berichte Verhandl. Akad. Leipzig {\bf 88}, 119--132 (1936)

\refno\CarminatiIKS.
Carminati, L., Iochum, B., Kastler, D., Sch\"ucker, T.:
On Connes' new gauge principle of general relativity. Can spinors
hear the forces of spacetime?
Preprint hep-th/9612228

\refno\ChamConnes.
Chamseddine, A. H., Connes, A.:
Universal formula for noncommutative geometry actions:
Unification of gravity and the Standard Model.
Phys. Rev. Lett. {\bf 77}, 4868--4871 (1996)

\refno\ConnesA.
Connes, A.:
The action functional in noncommutative geometry.
Commun. Math. Phys. {\bf 117}, 673--683 (1988)

\refno\ConnesNCGR.
Connes, A.:
Noncommutative geometry and reality.
J. Math. Phys. {\bf 36}, 6194--6231 (1995)

\refno\ConnesGrav.
Connes, A.:
Gravity coupled with matter and the foundation of noncommutative
geometry.
Preprint hep-th/9603053

\refno\ElizaldeVZ.
Elizalde, E., Vanzo, V., Zerbini, S.:
Zeta-function regularization, the multiplicative anomaly and the
Wodzicki residue.
Preprint hep-th/9701160

\refno\EstradaC.
Estrada, R.:
The Ces\`aro behaviour of distributions.
Preprint: San Jos\'e 1996

\refno\EstradaF.
Estrada, R., Fulling, S. A.:
The asymptotic expansion of spectral functions.
Pre\-print: College Station 1997

\refno\EstradaKfp.
Estrada, R., Kanwal, R. P.:
Regularization, pseudofunction and Hadamard finite part.
J. Math. Anal. Appl. {\bf 141}, 195--207 (1989)

\refno\EstradaKbook.
Estrada, R., Kanwal, R. P.:
Asymptotic Analysis: a Distributional Approach.
Bos\-ton: Birkh\"auser 1994

\refno\Amalthea.
Figueroa, H.:
Function algebras under the twisted product.
Bol. Soc. Paran. Mat. {\bf 11}, 115--129 (1990)

\refno\FullingA.
Fulling, S. A.:
The local geometric asymptotics of continuum eigenfunction expansions.
I. Overview.
SIAM J. Math. Anal. {\bf 13}, 891--912 (1982)

\refno\GelfandShilov.
Gelfand, I. M., Shilov, G. E.:
Generalized Functions I.
New York: Academic Press 1964

\refno\Gilkey.
Gilkey, P. B.:
Invariance Theory, the Heat Equation and the Atiyah--Singer Theorem,
2nd edition.
Boca Raton: CRC Press 1995

\refno\GrossmannLS.
Grossmann, A., Loupias, G., Stein, E. M.:
An algebra of pseudodifferential operators and quantum mechanics
in phase space.
Ann. Inst. Fourier (Grenoble) {\bf 18}, 343--368 (1968)

\refno\Grubb.
Grubb, G.:
Functional Calculus of Pseudodifferential Boundary Problems.
Boston: Birkh\"auser 1986

\refno\Gurarie.
Gurarie, D.:
The inverse spectral problem.
{\it in} Forty More Years of Ramifications: Spectral Asymptotics and
its Applications, Fulling, S. A., Narcowich, F. J., eds.,
Texas A\&M University, College Station (1991), pp.~77--99

\refno\Hardy.
Hardy, G. H.:
Divergent Series.
Oxford: Clarendon Press 1949

\refno\HormanderS.
H\"ormander, L.:
The spectral function of an elliptic operator.
Acta Math. {\bf 121}, 193--218 (1968)

\refno\HormanderIII.
H\"ormander, L.:
The Analysis of Linear Partial Differential Operators III.
Berlin: Springer 1985

\refno\IochumKS.
Iochum, B., Kastler, D., Sch\"ucker, T.:
On the universal Chamseddine--Connes action.
I. Details of the action computation.
Preprint hep-th/9607158

\refno\KalauW.
Kalau, W., Walze, M.:
Gravity, noncommutative geometry and the Wodzicki resi\-due.
J. Geom. Phys. {\bf 16}, 327--344 (1995)

\refno\KastlerEH.
Kastler, D.:
The Dirac operator and gravitation.
Commun. Math. Phys. {\bf 166}, 633--643 (1995)

\refno\Lagrange.
Lagrange, J.-L.:
Nouvelle m\'ethode pour r\'esoudre les \'equations litt\'erales
par la moyen des s\'eries.
M\'em. Acad. Royale des Sciences et
Belles-lettres de Berlin {\bf 24}, 251--326 (1770)

\refno\Gianni.
Landi, G.:
An introduction to noncommutative spaces and their geometry.
Preprint hep-th/9701078

\refno\Lojasiewicz.
{\Lpolish}ojasiewicz, S.:
Sur le valeur et la limite d'une distribution en un point.
Studia Math. {\bf 16}, 1--36 (1957)

\refno\Cordelia.
Mart\'{\i}n, C. P., Gracia-Bond\'{\i}a, J. M., V\'arilly, J. C.:
The Standard Model as a noncommutative geometry: the low energy
regime.
Preprint hep-th/9605001, Phys. Rep., to appear

\refno\McKeanS.
McKean, H. P., Singer, I. M.:
Curvature and the eigenvalues of the Laplacian.
J. Diff. Geom. {\bf 1}, 43--69 (1967)

\refno\Mulholland.
Mulholland, H. P.:
An asymptotic expansion for
$\sum_0^\infty (2n + 1) e^{-\sigma(n+\shalf)^2}$.
Proc. Camb. Philos. Soc. {\bf 24}, 280--289 (1928)

\refno\Rosenberg.
Rosenberg, S.:
The Laplacian on a Riemannian manifold.
Cambridge: Cambridge University Press 1997

\refno\Solomyak.
Solomyak, M. Z.:
Asymptotics of the spectrum of the Schr\"odinger operator with
nonregular homogeneous potential.
Math. USSR Sbornik {\bf 55}, 19--37 (1986)

\refno\Sirius.
V\'arilly, J. C., Gracia-Bond\'{\i}a, J. M.:
Connes' noncommutative differential geometry and the Standard Model.
J. Geom. Phys. {\bf 12}, 223--301 (1993)

\refno\Wodzicki.
Wodzicki, M.:
Noncommutative residue I: Fundamentals.
In Manin, Yu.~I. (ed.) $K$-theory, Arithmetic and Geometry.
Lecture Notes in Mathematics Vol.~1289, pp.~320--399.
Berlin: Springer 1987

\bye